\documentclass[11pt]{amsart}

\usepackage{amsmath,amssymb,latexsym,amsthm,enumerate}
\usepackage{verbatim,color}
\usepackage{comment}
\usepackage{hyperref}

\usepackage{graphicx}

\theoremstyle{plain}
\newtheorem{theorem}{Theorem}
\newtheorem{lem}[theorem]{Lemma}
\newtheorem{cor}[theorem]{Corollary}

\theoremstyle{remark}
\newtheorem*{remark}{Remark}
\newtheorem*{remarks}{Remarks}

\def\BEN{\begin{enumerate}}  \def\BI{\begin{itemize}}
\def\EEN{\end{enumerate}}   \def\EI{\end{itemize}}
   \def\sec{\section} \def\nn{\nonumber}

\def\beq{\begin{eqnarray}} \def\eeq{\end{eqnarray}}

\def\eqn#1{\begin{equation}#1\end{equation}}

\def\al*#1{\begin{align*}#1\end{align*}}

\def\ga*#1{\begin{gather*}#1\end{gather*}}

\def\alat*#1#2{\begin{alignat*}{#1}#2\end{alignat*}}
\def\bea{\begin{eqnarray*}}
\def\eea{\end{eqnarray*}}
\def\ml*#1{\begin{multline*}#1\end{multline*}}

\def\mbb{\mathbb} \def\mbf{\mathbf} \def\mrm{\mathrm}
\def\mc{\mathcal} \def\unl{\underline} \def\ovl{\overline}

\def\P{{\mathbb P}} \def\le{\left} \def\ri{\right} \def\i{\infty}
\def\C{{\mathbb C}}
\def\E{{\mathbb E}}   \def\R{{\mathbb R}}

 \def\EE{\mathcal{E}} 
\def\FF{\mathcal{F}}

\def\te#1{\mathrm{e}^{#1}}   

\def\WH{\widehat}

\def\I{\int}     \def\a{\alpha} \def\b{\beta}
\def\g{\gamma}  \def\d{\delta}   \def\th{\theta}

  \def\nn{\nonumber}   \def\s{\sigma}
\def\t{\tau}     \def\ps{\psi}
  \def\q{\qquad} 
\def\F{\Phi} \def\G{\Gamma}  \def\O{\Omega} 
  \def\td{\text{\rm d}}
\numberwithin{equation}{section}

\newcommand{\exit}{{\mbox{\, \vspace{3mm}}} \hfill\mbox{$\square$}}
\begin{document}

\title[On the drawdown of completely asymmetric L\'evy processes]
{On the drawdown of completely asymmetric L\'evy processes}
\thanks{{\it Acknowledgements:} We thank two anonymous referees for their useful comments
that led to improvements of the paper. MRP was supported in part by EPSRC
Mathematics Platform Grant EP/I019111/1 and the NWO-STAR cluster.}

\author{Aleksandar Mijatovi\'{c}}
\address{Department of Statistics, University of Warwick, UK}
\email{a.mijatovic@warwick.ac.uk}

\author{Martijn R. Pistorius}
\address{Department of Mathematics, Imperial College London, UK \newline
and Korteweg-de Vries Institute for Mathematics,
University of Amsterdam} \email{m.pistorius@imperial.ac.uk}



\begin{abstract}
The {\em drawdown} process $Y$ of a completely asymmetric L\'{e}vy
process $X$ is equal to $X$ reflected at its running supremum
$\overline{X}$: $Y = \overline{X} - X$. In this paper
we explicitly express in terms of the scale function and the L\'{e}vy measure of $X$
the law of the sextuple of the first-passage time of $Y$ over the level $a>0$,
the time $\overline{G}_{\tau_a}$ of the last supremum of $X$
prior to $\tau_a$, the infimum $\unl X_{\tau_a}$ and supremum $\ovl X_{\tau_a}$
of $X$ at $\tau_a$ and the undershoot $a - Y_{\tau_a-}$ and overshoot
$Y_{\tau_a}-a$ of $Y$ at $\tau_a$. As application we obtain
explicit expressions for the laws of a number of functionals
of drawdowns and rallies in a  completely asymmetric exponential L\'{e}vy model.
\end{abstract}

\maketitle

\sec{Introduction}\label{sec:intro} A completely asymmetric L\'evy
process is a real-valued stochastic process with c\`adl\`ag paths
that has independent stationary increments whose jump sizes all
have the same sign. Its {\it drawdown process}, also known as the
{reflected} process, is the difference of its running supremum and
its current value. Closely related is the {\it rally} process
which is defined as the difference of the current value and the
running infimum, and is equal to the drawdown process of the
negative of the process. 
This paper is concerned with a {\em distributional} study of the
drawdown process of a completely asymmetric L\'{e}vy process 
For such a
L\'{e}vy process~$X$ the law is identified of the following sextuple concerning
its drawdown process $Y$: \begin{equation}(\t_a,\ovl
G_{\tau_a}, \ovl X_{\t_a}, \unl X_{\t_a}, a- Y_{\t_a-}, Y_{\t_a}-a),
\q\text{where}\q Y_t = \sup_{0\leq s\leq t} X_s - X_t, \quad t\ge
0.\end{equation} Here the components of the vector are given by
$\tau_a$, the first-passage time of $Y$ over a level $a>0$, $\ovl
G_{\tau_a}$, the last time that $X$ is at its supremum prior to
$\tau_a$, $\ovl X_{\tau_a}$ and~$\unl X_{\tau_a}$, the running supremum and infimum
at $\tau_a$, and $a- Y_{\t_a-}$ and $Y_{\t_a}-a$, the undershoot and the overshoot of
$Y$ at the epoch $\tau_a$.

The drawdown process has been the object
of considerable interest in various areas of applied probability.
It has been studied for instance in queueing theory ({\it e.g.}, Asmussen
\cite{APQ}), risk theory and mathematical genetics.
The drawdown process has also been employed in financial
modelling, in the construction of tractable,
path-dependent risk/performance measures. In the
context of real estate portfolio optimisation, Hamelink \&
Hoesli~\cite{HH} considered the running maximum of the drawdown
process as an investment performance criterion. Checkhlov et
al.~\cite{CUZ} introduced a one-parameter family of portfolio risk measures that was
called conditional drawdown and defined to be equal to the mean of
a percentage of the worst portfolio drawdowns.
Pospisil et al.~\cite{PosVec} proposed the probability of a
drawdown of a given size occurring before a rally of a given size
as risk measure, and calculated this probability in the setting of
one-dimensional diffusion models; the finite horizon case for
Brownian motion was treated by Zhang \& Hadjiliadis~\cite{ZH}.

A number of papers has been devoted to a distributional study of
functionals of the drawdown process. The joint Laplace of the time
to a given drawdown and the running maximum of a Brownian motion
with drift was derived by Taylor~\cite{Taylor}; this joint law was
obtained  by Lehoczky~\cite{Lehoczky} in the case of a general diffusion. An
explicit expression for the expectation and the density of the
maximum drawdown of Brownian motion was derived by Douady~{\it et
al.}~\cite{DSY}; the case of Brownian motion with drift was covered
by Magdon~{\it et al.}~\cite{MAPA} where also the large time
asymptotics of the expectation were derived.

The drawdown process also features in the solution of a number of optimal
investment problems. Under the geometric Brownian motion model
the optimal  time to exercise the Russian
option, which pays out the largest historical value of the stock at the
moment of exercise,
was shown by Shepp \& Shiryaev \cite{SHSH} to be
given by the first-passage of a drawdown
process over a certain constant level. Such a first-passage time
is also optimal when linear cost is included
(Meilijson~\cite{Meilijson}), or under a spectrally negative L\'{e}vy model for the stock price
(Avram {\it et al.}~\cite{AKP}).

Although still widely used as benchmark, mainly on account of its
analytical tractability, it is by now well established that many features
of Samuelson's classical geometric Brownian motion model for
the price of a stock are not supported by empirical data.
A class of tractable models that captures typical features of
stock returns data such as fat tails, asymmetry and excess
kurtosis is that of exponential L\'{e}vy processes. This class has
received considerable attention in the literature
---we refer to Cont \& Tankov~\cite{ContTankov}
and Boyarchenko \& Levendorskii~\cite{BoL}
for background and references.
By restricting ourselves to L\'{e}vy processes with jumps of a
single sign, we are able to draw on the fluctuation theory for this class of stochastic processes, which is considerably more explicit than
in the case of general L\'{e}vy processes. Empirical support for a
model from this class was given in Carr~\&~Wu~\cite{CarrWu}, where options on
the S\&P 500 index were studied. Carr~\&~Wu~\cite{CarrWu}
demonstrated that the finite-moment log-stable model, which is an
exponential L\'{e}vy model driven by a spectrally negative stable
process, provided a good fit to quoted option prices across
maturities.

By way of application, we employ the sextuple law to obtain semi-analytical expressions for the expectations of a number of path-functionals of the drawdown process of an exponential L\'{e}vy process, which provide a description of different aspects of the riskiness of the model: (i) the probability that, on a given time horizon, a new minimum is attained (e.g. by a jump) at the first moment of a drawdown of a given size;
(ii) the expected size of the drawdown process at the first moment that a drawdown of a given size occurs, given that this happens before a finite time-horizon and (iii) the probability that, on a finite time-horizon, a drawdown of a given size occurs before a rally of a given size.
These and related explicit expressions can form the basis for the analysis of optimal investment problems involving the drawdown process---in the interest of brevity,
such investigations are left for future research.

In the literature two approaches have been successfully adopted to
identify the distributions of path-functionals of the
drawdown process: the {\it martingale} approach, exemplified in
e.g.~\cite{CKO,NY}, which rests on the identification of certain
martingales, and the {\it excursion-theoretic} approach,
employed  in a.o.~\cite{AKP,Baurdoux09,Chamont_Doney,Doney,KPR,Pe},
which is based on the theorem by It\^o stating that the process of
excursions of $Y$ away from zero forms a Poisson point process. In
this paper we will follow the latter approach.

The remainder of the paper is organized as follows. In
Section~\ref{sec:prel} preliminaries are reviewed and the notation
is set. In Sections~\ref{subsec:quintuple} and~\ref{subsec:triple}
the law of the sextuple is derived for a spectrally negative and
spectrally positive L\'{e}vy process, respectively.
In Section~\ref{sec:DD} the sextuple law is employed to derive
analytically explicit identities for a number of characteristics of drawdowns
and rallies in an exponential L\'{e}vy model.  The proofs of the main
theorems are contained in Section~\ref{sec:proofs}.

\sec{Preliminaries}\label{sec:prel} In this section we set the
notation. For the background on the fluctuation theory of
spectrally negative L\'{e}vy processes refer to
Bertoin \cite[Chapter~VII]{bert96} and Kyprianou \cite{Kyprianou}.

Let $X=(X_t)_{t\ge0}$ be a L\'{e}vy process defined on
$(\O,\FF,(\FF_t)_{t\ge0},\P)$, a filtered probability space which
satisfies the usual conditions.  Assume throughout the paper that
$X$ is spectrally negative. Under $\P$ we have $X_0=0$ $\P$-a.s.
For any $x\in\R$ we denote by $\P_x$ the law of the L\'evy process
started from $x$, i.e. the law of the process $x+X$ under $\P$. To
avoid the case of trivial reflected processes we exclude $X$ that
has monotone paths, i.e. $X$ is assumed to be neither the negative
of a subordinator nor a deterministic drift upwards. Since the
jumps of $X$ are all non-positive, the moment generating function
$\E[\te{\th X_t}]$ exists for all $\th\ge0$ and is given by
$\ps(\th) = t^{-1} \log \E[\te{\th X_t}]$ for some function
$\ps(\th)$. The function $\psi$ which is well defined at least on
the positive half-axis where it is strictly convex with the
property that $\lim_{\th\to\i}\ps(\th)=+\i$. Let $\F(0)$ be the
largest root of $\ps(\th)=0$. On $[\F(0),\i)$ the function $\ps$
is strictly increasing and we denote its right-inverse function by
$\F:[0,\i)\to[\F(0),\i)$.

For $q\ge0$, there exists a continuous increasing function
$W^{(q)}: [0,\i) \to [0,\i)$, called the {\it $q$-scale function},
with Laplace transform \eqn{\label{eq:defW} \I_0^\i \te{-\th x}
W^{(q)} (x) \td x = (\ps(\th) - q)^{-1},\q\q\th > \F(q). } The
function $W^{(q)}$ is extended to $x\in(-\infty,0)$ by
$W^{(q)}(x)=0$. A related $q$-scale function $Z^{(q)}$ is defined
by
\begin{equation}
\label{eq:defZ} Z^{(q)}(x) = 1 + q \I_0^xW^{(q)}(z)\td
z,\qquad\text{for $x\in\R$.}
\end{equation}
The $q$-scale function $W^{(q)}$ is left- and right-differentiable
on $(0,\i)$ and we denote the right- and left-derivative of
$W^{(q)}$ by $W^{(q)\prime}_+$ and $W^{(q)\prime}_-$ respectively
for any $q\geq0$. Furthermore if a Gaussian component is present
then for any $q\geq0$ we have $W^{(q)}\in C^2(0,\i)$
(see~\cite{CK}). In this case we denote by $W^{(q)\prime}$ and
$W^{(q)\prime\prime}$ the first and the second derivative of the
scale function $W^{(q)}$ respectively.

We next briefly review two-sided exit results which will be employed in the sequel.
Let $u<v$ and $x\in[u,v]$ and define the
first-passage times $T^-_u$ and $T^+_v$ as
$$
T^-_u = \inf\{t\ge 0: X_t < u\}
\quad\text{ and }\quad
T^+_v = \inf\{t\ge 0: X_t > v\}.
$$
Let
$T_{u,v} = T^-_u\wedge T^+_v$
be the first time
$X$
started at
$x$
enters the set
$\R\setminus[u,v]$. The two sided first-passage results
\cite{bert97} read
\begin{eqnarray}
\label{eq:two-sided1} \E_x\le[\te{-q T_{u,v}}I_{\le\{X_{T_{u,v} = v}\ri\}}\ri] &=&
\frac{W^{(q)}(x-u)}{W^{(q)}(v-u)},\q x\in[u,v],\\
\label{eq:two-sided2}
\E_x\le[\te{-q T_{u,v}}I_{\le\{X_{T_{u,v} \leq u}\ri\}}\ri] &=&
Z^{(q)}(x-u) - Z^{(q)}(v-u)\frac{W^{(q)}(x-u)}{W^{(q)}(v-u)},\q x\in[u,v].
\end{eqnarray}
%

A L\'{e}vy process started at zero \textit{creeps downwards}
(resp. \textit{upwards}) over a level $x<0$ (resp. $x>0$) if, with
positive probability, the first time that it enters the set
$(-\infty,x)$ (resp. $(x,\infty)$) this does not happen by a jump.
It is an immediate consequence of the Wiener-Hopf factorization
that a L\'{e}vy process creeps both upwards and downwards if and
only if its Gaussian component is not zero. Therefore a spectrally
negative L\'{e}vy process with non-monotone trajectories creeps
downwards if and only if it has a positive Gaussian coefficient
(e.g. \cite[p. 175]{bert96}). In this case the results of
Millar~\cite{Millar} imply that
\begin{equation}
\E_x\le[\te{-q T^-_u}I_{\{X_{T^-_u}=u\}}\ri] = \frac{\s^2}{2}\le(W^{(q)\prime}(x-u) - \F(q)W^{(q)}(x-u)\ri)\quad
\text{for any}\quad u\leq x,\,\, q\geq0.
\label{eq:creep}
\end{equation}
Both sides of the equality in~\eqref{eq:creep} are understood to
be equal to zero if $\s=0$. The formula for the probability that
$X$ leaves the interval $[u,v]$ by hitting $u$ follows
from~\eqref{eq:two-sided1},~\eqref{eq:creep} and the strong Markov
property and is given by the following expression:
\begin{equation}
\label{eq:twosmoothdown}
\E_{x}\le[\te{-qT_{u,v}}I_{\left\{X_{T_{u,v}}=u\right\}}\ri] = \frac{\s^2}{2}\le(
W^{(q)\prime}(x-u) - \frac{W^{(q)\prime}(v-u)}{W^{(q)}(v-u)}W^{(q)}(x-u)\ri)
\end{equation}
for
$x\in[u,v]$.
Again the expression is understood to be equal to $0$ if $\s=0$.

\section{Reflected spectrally negative L\'{e}vy processes}
\label{sec:exc}
Let
$\ovl X_t = \sup_{0\leq u\leq t}\le\{X_u\ri\}$
and
$\unl X_t = \inf_{0\leq u\leq t}\le\{X_u\ri\}$
and define the reflected processes
$Y=(Y_t)_{t\geq0}$
and
$\WH Y=(\WH Y_t)_{t\geq0}$
by
\begin{eqnarray}
\label{eq:Def_Reflected_Procs}
Y_t = \ovl X_t - X_t &\quad\text{and}\quad & \WH Y_t =X_t- \unl X_t.
\end{eqnarray}
The focus of this paper are the first-passage times over a level
$a>0$ of the reflected processes $Y$ and $\WH Y$,
\begin{eqnarray}
\label{eq:Def_First_Passige_a} \tau_a := \inf\{t\geq 0: Y_t > a\}
\qquad\text{and}\qquad \WH \tau_a := \inf\{t\geq 0: \WH Y_t > a\}.
\end{eqnarray}
It is well known that the stopping times $\tau_a$ and $\WH \tau_a$
are finite $\P$-a.s. Furthermore, in \cite[Theorem 1]{AKP} the
Laplace transform of $\tau_a$ was identified as
\begin{equation}\label{eq:taua} \E[\te{-q\tau_a}] = Z^{(q)}(a) - q
\frac{W^{(q)}(a)^2}{W^{(q)\prime}_+(a)}.
\end{equation}
The last times before time $t$ that $X$ visits its running
supremum and infimum are  denoted by $\ovl G_{t}$ and $\unl G_{t}$
respectively, where
\begin{equation}
\label{eq:Last_time_Def}
\ovl G_{t}  =  \sup \left\{s\leq t\,:\, X_{s}\text{ or } X_{s-}=\ovl X_{s}\right\}\quad
\text{and}\quad
\unl G_t = \sup\{s \leq t: X_s\text{ or } X_{s-}=\unl X_s\}.
\end{equation}
Note that
$\ovl G_{t}$
(resp.
$\unl G_{t}$)
can be viewed as the last time before time $t$
that the reflected process
$Y$
(resp.
$\WH Y$) is equal to
$0$.

In Section~\ref{subsec:quintuple} we characterise the joint law of the following sextuple
of random variables:
$$\begin{array}{lll}
\tau_a & &  \text{the first-passage time over a level $a$ of the reflected process $Y$,}\\
\ovl G_{\tau_a} & &
\text{the last time that $X$ was at its
running supremum
prior to the first-passage time $\tau_a$,}\\
\ovl X_{\tau_a} & &
\text{the
supremum of
$X$
at the first-passage time $\tau_a$,}\\
\unl X_{\tau_a} & &
\text{the infimum  of $X$ at the first-passage time $\tau_a$,}\\
Y_{\tau_a-} & & \text{the position of the reflected process just before it crosses the
level $a$,}\\
Y_{\tau_a}-a & & \text{the overshoot of the reflected process $Y$
over the level $a$.}
\end{array}$$
In
Section~\ref{subsec:triple} we give the joint law of the following quadruple of
random variables
$$\begin{array}{lll}
\WH \tau_a & &  \text{the first-passage time over a level $a$  of the reflected process $\WH Y$,}\\
\unl G_{\WH \tau_a} & &
\text{the last time that
$X$ was at its running infimum
prior to the first-passage time $\WH\tau_a$},\\
\ovl X_{\WH \tau_a} & &
\text{the
supremum of
$X$
at the first-passage time $\WH\tau_a$,}\\
\unl X_{\WH \tau_a} & & \text{the infimum  of $X$ at the
first-passage time $\WH\tau_a$.}
\end{array}$$
Note that in this case, since $X$ is assumed to be spectrally
negative, the reflected process $\WH Y$ can only jump down. Since
$\WH Y$ is right-continuous with left limits at the first-passage
time $\WH \tau_a$ we have $\WH Y_{\WH \tau_a-}=\WH Y_{\WH
\tau_a}=a$ a.s.

\subsection{The sextuple law}
\label{subsec:quintuple} We now give the law of the sextuple
$\left(\tau_a, \ovl G_{\tau_a},\ovl X_{\tau_a},\unl X_{\tau_a},
Y_{\tau_a-}, Y_{\tau_a}-a\right)$. Define for any $a>0$ and
$p,q\geq0$ the map $F_{p,q,a}:\mbb R_+\to\mbb R_+$ by
\begin{eqnarray}
\label{Def:F} F_{p,q,a}(y) =
\lambda(a,q) \exp\le(-
y\lambda(a,p)\ri), \qquad y\in\mbb R_+,
\end{eqnarray}
where $\lambda(a,q)$ is the ratio of the derivative of the $q$-scale function
and the $q$-scale function at $a$,
\begin{equation}\label{eq:laq}
\lambda(a,q) = \frac{W^{(q)\prime}_+(a)}{W^{(q)}(a)}.
\end{equation}

 Let $\Lambda$ be the L\'{e}vy measure and $\s^2\geq0$
the Gaussian coefficient of $X$. Denote by 
$R^{(q)}_a(\td y)= \E\left[\int_0^{\tau_a}\te{-qt}I_{\{Y_t\in\td y\}}\td t\right]$ 
the $q$-resolvent measure of $Y$ killed upon first exit from $[0,a]$,
which may be expressed in terms of the $q$-scale function $W^{(q)}$
by (\cite[Theorem 1]{P})
\begin{eqnarray}\label{eq:res}
R^{(q)}_a(\td y) = \left[ \lambda(a,q)^{-1}\,
W^{(q)}(\td y) - W^{(q)}(y)\td y\right], \qquad y\in[0,a],
\end{eqnarray}
Furthermore, let $\Delta^{(q)}(a)$ denote the
expression that is given in terms of the derivatives of $W^{(q)}$ by
\begin{equation}
\Delta^{(q)}(a) = \frac{\sigma^2}{2}
\left[W^{(q)\prime}(a) -
\lambda(a,q)^{-1} W^{(q)\prime\prime}(a)\right], \label{eq:Daq}
\end{equation}
where the expression on the right-hand side of the
equality in~\eqref{eq:Daq} is taken to be zero if
$\s=0$. As we will see below, the Laplace transform
of $\tau_a$ on the event that $Y$ creeps over the level $a$
is equal to $\Delta^{(q)}(a)=\E[\te{-q\tau_a}I_{\{Y_{\tau_a}=a\}}]$.

\begin{theorem}\label{thm:lawreflect1}
Let $X_0=x\in\R$ and $a>0$
and define events
$$A_o=\left\{ \unl
X_{\tau_a}\geq u, \ovl X_{\tau_a}\in\td v, Y_{\tau_a-}\in\td y,
Y_{\tau_a}-a\in\td h \right\}\quad\text{and}\quad A_c=\left\{ \unl
X_{\tau_a}\geq u, \ovl X_{\tau_a}\in\td v, Y_{\tau_a}=a\right\},$$
where $u,v,y$ and $h$ satisfy
\begin{equation}
\label{eq:ParmRest} u\leq x,\> y\in[0,a],\> v\geq x\vee (u+a)
\quad\text{and}\q h\in(0,v-u-a].
\end{equation}
Then for any $q,r\geq0$ the following identities hold true:
\begin{eqnarray}
\nonumber
\lefteqn{\E_x\left[\te{-q\t_a-r\ovl G_{\tau_a}}I_{A_o}\right] \> =} \\
&  & \label{eq:QuitLaq_overshoot} \frac{W^{(q+r)}((x-u)\wedge
a)}{W^{(q+r)}(a)}\,F_{q+r,q,a}\left(v-(x\vee(u+a))\right)\,R_a^{(q)}(\td
y)\,\Lambda\left(y-a-\td h\right),
\\
\nonumber
 \lefteqn{\E_x\left[\te{-q\t_a-r\ovl G_{\tau_a}}I_{A_c}\right]\> =} \\
&& \frac{W^{(q+r)}((x-u)\wedge a)}{W^{(q+r)}(a)}\,
F_{q+r,q,a}\left(v-(x\vee(u+a))\right)  \Delta^{(q)}(a)
\label{eq:QuitLaq_Creeping}
\end{eqnarray}
where $I_{\{\cdot\}}$ denotes the indicator of the set
${\{\cdot\}}$ and $c\wedge d=\min\{c,d\}$, $c\vee d=\max\{c,d\}$ for 
$c,d\in\R$.
\end{theorem}

\begin{remarks}
\noindent (i)  The
formulas in~\eqref{eq:QuitLaq_overshoot}
and~\eqref{eq:QuitLaq_Creeping} determine the law of the sextuple
$\left(\tau_a, \ovl G_{\tau_a},\ovl X_{\tau_a},\unl X_{\tau_a},
Y_{\tau_a-}, Y_{\tau_a}-a\right)$.
Indeed, note that if the parameters $u,v,y$ and $h$ do not
satisfy the restrictions in~\eqref{eq:ParmRest}, then we have
$\P_x\left[A_o\right]=\P_x\left[A_c\right]=0$.  In particular, if we take
$u=-\infty$, then~\eqref{eq:ParmRest} places a restriction neither
on the size $h\in(0,\infty)$ of the overshoot of the reflected
process nor on the level $v\geq x$ of the supremum of $X$ attained
at $\tau_a$. This is also transparent from Figure~\ref{fig:Main}.
\smallskip

\noindent (ii) Note that on the event
$A_o$
we have
$X_{\tau_a}=v-a-h$ $\P_x$-a.s.
and hence
$\unl X_{\tau_a}=(v-a-h)\wedge\unl X_{\tau_a-}$.
Define the event
$$A_o^-=\left\{ \unl
X_{\tau_a-}\geq u, \ovl X_{\tau_a}\in \td v, Y_{\tau_a-}\in\td y,
Y_{\tau_a}-a\in\td h \right\}.
$$
Then the following implications hold:
\begin{eqnarray*}
h\in(0,v-u-a] & \Longrightarrow  &
A_o^-=A_o\quad
\P_x\text{-a.s.}\\
 \quad
h>v-u-a & \Longrightarrow   &
\unl X_{\tau_a} = v-a-h\quad
\P_x\text{-a.s. on } A_o^-\quad\text{and}\quad
\P_x\left[A_o\right]=0.
\end{eqnarray*}
Furthermore, it is clear from step~(2) of the  proof
of Theorem~\ref{thm:lawreflect1} (presented in Section~\ref{sec:proofs})
that
$\E_x\left[\te{-q\t_a-r\ovl G_{\tau_a}}I_{A_o^-}\right]$
is given by the formula in~\eqref{eq:QuitLaq_overshoot},
which in this case holds without a restriction on the size
of the overshoot
$h$
from~\eqref{eq:ParmRest}.
\smallskip

\noindent (iii) Note that the reflected process $Y$ {creeps} over
the level $a$ (i.e. $\P[Y_{\tau_a}=a]>0$) if and only if $X$
creeps downwards which is the case precisely if its Gaussian
coefficient $\sigma^2$ is strictly positive. 
Hence the expression
$\Delta^{(q)}(a)$
given in~\eqref{eq:Daq} is well-defined, as
$W^{(q)}$ is twice differentiable at $a>0$ if $\sigma^2$ is
strictly positive.
\smallskip

\noindent (iv) Schematic representations of a typical path
of the process
$X$
in the events
$A_o$ and $A_c$
of Theorem~\ref{thm:lawreflect1}
are given in the bottom and top pictures in
Figure~\ref{fig:Main}.
Note that
$A_c$
contains the paths
of the reflected process
$Y$
when it creeps over the level
$a$
and that
$Y_{\tau_a-}=a$
$\P$-a.s. on
$A_c$.
Similarly
$A_o$
consists of the paths
of
$Y$
which first enter
$(a,\infty)$
by a jump.

\smallskip
\begin{figure}[p]
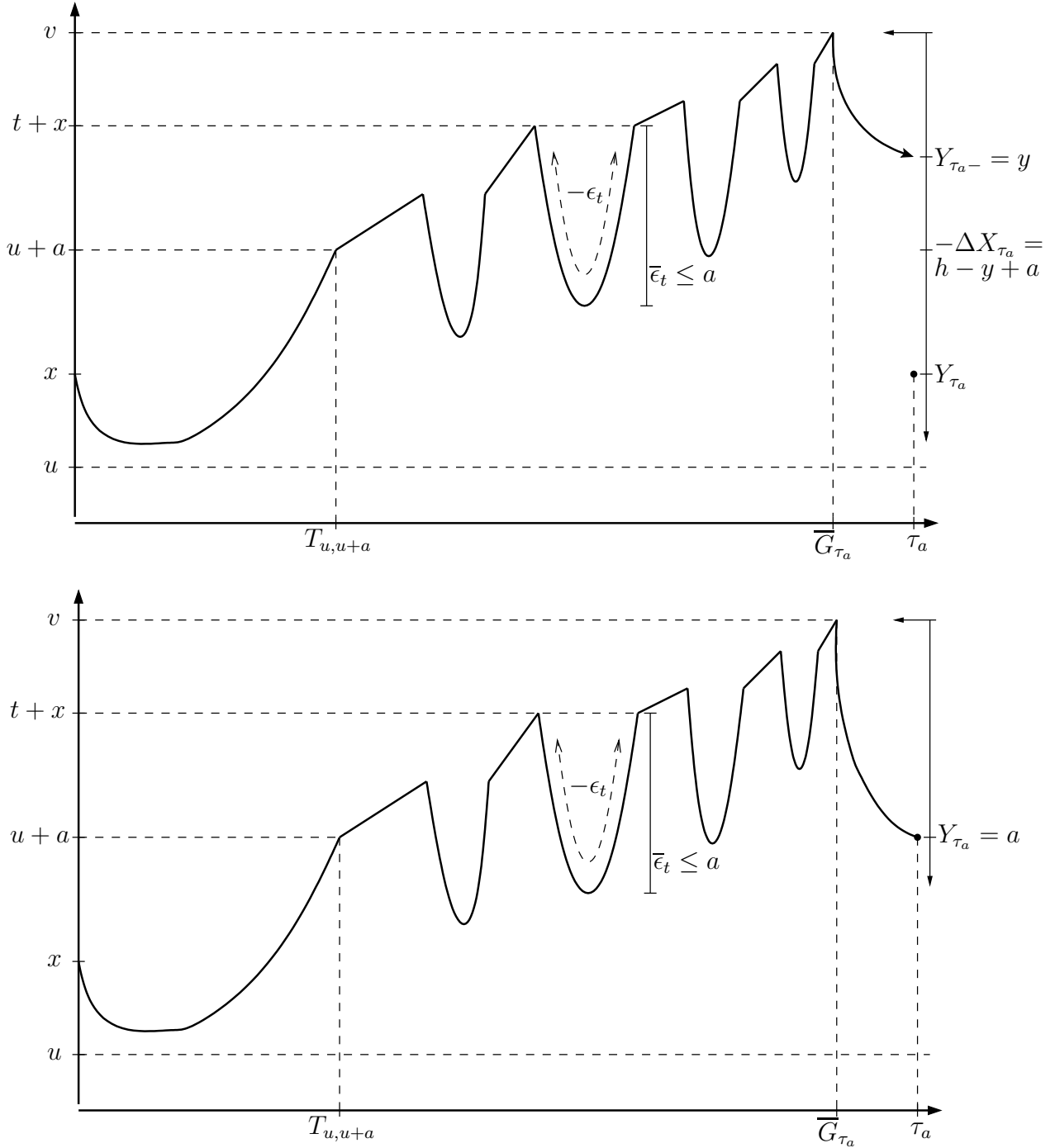

\begin{picture}(0,0)%
\includegraphics{Path.pstex}%
\end{picture}%
\setlength{\unitlength}{4144sp}%
\begingroup\makeatletter\ifx\SetFigFont\undefined%
\gdef\SetFigFont#1#2#3#4#5{%
  \reset@font\fontsize{#1}{#2pt}%
  \fontfamily{#3}\fontseries{#4}\fontshape{#5}%
  \selectfont}%
\fi\endgroup%
\begin{picture}(7507,4105)(616,-3694)
\put(901,-2356){\makebox(0,0)[lb]{\smash{{\SetFigFont{12}{14.4}{\rmdefault}{\mddefault}{\updefault}{\color[rgb]{0,0,0}$x$}%
}}}}
\put(901,-3031){\makebox(0,0)[lb]{\smash{{\SetFigFont{12}{14.4}{\rmdefault}{\mddefault}{\updefault}{\color[rgb]{0,0,0}$u$}%
}}}}
\put(631,-1456){\makebox(0,0)[lb]{\smash{{\SetFigFont{12}{14.4}{\rmdefault}{\mddefault}{\updefault}{\color[rgb]{0,0,0}$u+a$}%
}}}}
\put(901,119){\makebox(0,0)[lb]{\smash{{\SetFigFont{12}{14.4}{\rmdefault}{\mddefault}{\updefault}{\color[rgb]{0,0,0}$v$}%
}}}}
\put(2791,-3571){\makebox(0,0)[lb]{\smash{{\SetFigFont{12}{14.4}{\rmdefault}{\mddefault}{\updefault}{\color[rgb]{0,0,0}$T_{u,u+a}$}%
}}}}
\put(5290,-1650){\makebox(0,0)[lb]{\smash{{\SetFigFont{12}{14.4}{\rmdefault}{\mddefault}{\updefault}{\color[rgb]{0,0,0}$\overline \epsilon_t\leq a$}%
}}}}
\put(7350,-2380){\makebox(0,0)[lb]{\smash{{\SetFigFont{12}{14.4}{\rmdefault}{\mddefault}{\updefault}{\color[rgb]{0,0,0}$Y_{\tau_a}$}%
}}}}
\put(7350,-811){\makebox(0,0)[lb]{\smash{{\SetFigFont{12}{14.4}{\rmdefault}{\mddefault}{\updefault}{\color[rgb]{0,0,0}$Y_{\tau_a-}=y$}%
}}}}
\put(7150,-3571){\makebox(0,0)[lb]{\smash{{\SetFigFont{12}{14.4}{\rmdefault}{\mddefault}{\updefault}{\color[rgb]{0,0,0}$\tau_a$}%
}}}}
\put(6481,-3616){\makebox(0,0)[lb]{\smash{{\SetFigFont{12}{14.4}{\rmdefault}{\mddefault}{\updefault}{\color[rgb]{0,0,0}$\overline G_{\tau_a}$}%
}}}}
\put(7350,-1606){\makebox(0,0)[lb]{\smash{{\SetFigFont{12}{14.4}{\rmdefault}{\mddefault}{\updefault}{\color[rgb]{0,0,0}$h-y+a$}%
}}}}
\put(7350,-1426){\makebox(0,0)[lb]{\smash{{\SetFigFont{12}{14.4}{\rmdefault}{\mddefault}{\updefault}{\color[rgb]{0,0,0}$-\Delta X_{\tau_a}=$}%
}}}}
\put(4681,-1051){\makebox(0,0)[lb]{\smash{{\SetFigFont{12}{14.4}{\rmdefault}{\mddefault}{\updefault}{\color[rgb]{0,0,0}$-\epsilon_t$}%
}}}}
\put(676,-556){\makebox(0,0)[lb]{\smash{{\SetFigFont{12}{14.4}{\rmdefault}{\mddefault}{\updefault}{\color[rgb]{0,0,0}$t+x$}%
}}}}
\end{picture}

\vspace{0.3cm}
\begin{picture}(0,0)%
\includegraphics{Path_2.pstex}%
\end{picture}%
\setlength{\unitlength}{4144sp}%
\begingroup\makeatletter\ifx\SetFigFont\undefined%
\gdef\SetFigFont#1#2#3#4#5{%
  \reset@font\fontsize{#1}{#2pt}%
  \fontfamily{#3}\fontseries{#4}\fontshape{#5}%
  \selectfont}%
\fi\endgroup%
\begin{picture}(7417,4105)(616,-3694)
\put(901,-2356){\makebox(0,0)[lb]{\smash{{\SetFigFont{12}{14.4}{\rmdefault}{\mddefault}{\updefault}{\color[rgb]{0,0,0}$x$}%
}}}}
\put(901,-3031){\makebox(0,0)[lb]{\smash{{\SetFigFont{12}{14.4}{\rmdefault}{\mddefault}{\updefault}{\color[rgb]{0,0,0}$u$}%
}}}}
\put(631,-1456){\makebox(0,0)[lb]{\smash{{\SetFigFont{12}{14.4}{\rmdefault}{\mddefault}{\updefault}{\color[rgb]{0,0,0}$u+a$}%
}}}}
\put(901,119){\makebox(0,0)[lb]{\smash{{\SetFigFont{12}{14.4}{\rmdefault}{\mddefault}{\updefault}{\color[rgb]{0,0,0}$v$}%
}}}}
\put(2791,-3571){\makebox(0,0)[lb]{\smash{{\SetFigFont{12}{14.4}{\rmdefault}{\mddefault}{\updefault}{\color[rgb]{0,0,0}$T_{u,u+a}$}%
}}}}
\put(5290,-1650){\makebox(0,0)[lb]{\smash{{\SetFigFont{12}{14.4}{\rmdefault}{\mddefault}{\updefault}{\color[rgb]{0,0,0}$\overline \epsilon_t\leq a$}%
}}}}
\put(7150,-3571){\makebox(0,0)[lb]{\smash{{\SetFigFont{12}{14.4}{\rmdefault}{\mddefault}{\updefault}{\color[rgb]{0,0,0}$\tau_a$}%
}}}}
\put(6481,-3616){\makebox(0,0)[lb]{\smash{{\SetFigFont{12}{14.4}{\rmdefault}{\mddefault}{\updefault}{\color[rgb]{0,0,0}$\overline G_{\tau_a}$}%
}}}}
\put(7350,-1461){\makebox(0,0)[lb]{\smash{{\SetFigFont{12}{14.4}{\rmdefault}{\mddefault}{\updefault}{\color[rgb]{0,0,0}$Y_{\tau_a}=a$}%
}}}}
\put(4681,-1096){\makebox(0,0)[lb]{\smash{{\SetFigFont{12}{14.4}{\rmdefault}{\mddefault}{\updefault}{\color[rgb]{0,0,0}$-\epsilon_t$}%
}}}}
\put(631,-556){\makebox(0,0)[lb]{\smash{{\SetFigFont{12}{14.4}{\rmdefault}{\mddefault}{\updefault}{\color[rgb]{0,0,0}$t+x$}%
}}}}
\end{picture}%
\caption{\footnotesize{The top figure contains a schematic representation of a typical path
of the process
$X$
in the event
$A_o=\left\{
\unl X_{\tau_a}\geq u,
\ovl X_{\tau_a}\in \td v,
Y_{\tau_a-}\in\td y,
Y_{\tau_a}-a\in\td h
\right\}$.
This figure illustrates the idea behind the proof of Theorem~\ref{thm:lawreflect1}:
the depicted path satisfies
$T_{u,u+a}=T_{u+a}^+$
and
$T_v^+=\ovl G_{\tau_a}$.
It is clear from the figure that the
trajectory can be decomposed into the following three parts,
each of which is analysed separately
in the proof of Theorem~\ref{thm:lawreflect1}:
(1) the segment over the interval
$[0,T_{u,u+a}]$,
(2) the jump at time
$\tau_a$
and
(3) the segment that straddles
$\ovl G_{\tau_a}$
over the time interval
$(T_{u,u+a},\tau_a)$.
The bottom figure depicts a schematic representation of a typical path
of the process
$X$
in the event
$A_c=\left\{
\unl X_{\tau_a}\geq u,
\ovl X_{\tau_a}\in \td v,
Y_{\tau_a}=a
\right\}$.
}}
\label{fig:Main}
\end{figure}
\end{remarks}

Explicit expressions, in terms
of the 0-scale function,
for the marginal distributions of $\ovl X_{\tau_a}$ and $\unl X_{\tau_a-}$
under 
$\mbb P_x$
can be obtained from Theorem~\ref{thm:lawreflect1}. 

\begin{cor}
\label{cor:Thm1} For any $a>0$ and $x\in\mbb R$, the following
hold true:
\begin{eqnarray}
\label{eq:Min_Law}
\P_x(x-\unl X_{\tau_a-} \leq z) & = & W(z\wedge a)/W(a), \qquad z\in\R_+,\\
\P_x(\ovl X_{\tau_a} - x\geq z) & = &  \exp\left(-z W_+'(a)/W(a)\right), \qquad z\in\R_+,
\label{eq:ExpLaw}
\end{eqnarray}
where $W=W^{(0)}$ denotes the $0$-scale function.
\end{cor}

\noindent {\it Proof }{\it of Corollary \ref{cor:Thm1}}
Setting $q$ and $r$ equal to $0$ and letting
$u=-\infty$ in~\eqref{eq:QuitLaq_overshoot}--\eqref{eq:QuitLaq_Creeping},
and integrating over $y\in[0,a]$, $h\in(0,\infty)$ and $v\in(x,\infty)$
the right-hand sides of the these expressions we obtain that $c(a)=1$ where
\begin{equation}\label{eq:ca}
c(a) := \int_0^\infty\int_0^a R_a^{(0)}(\td y)\Lambda(y-a-\td
h) + \Delta^{(0)}(a).
\end{equation}
Fixing $x,a, u_0, v_0$, such that $a>0$, $u_0+a\leq v_0$ and $x\in[u_0,v_0]$, and
again integrating the right-hand side
of~\eqref{eq:QuitLaq_overshoot} over
$y\in[0,a]$, $h\in(0,\infty)$ and $v\in(v_0,\infty)$ and the right-hand side of  Eqns.~\eqref{eq:QuitLaq_Creeping} over $v\ge
v_0$  we find
\begin{equation*}
\P_x(\unl X_{\tau_a-}\ge u_0, \ovl X_{\tau_a}\ge v_0) = c(a) \cdot
\frac{W^{(0)}((x-u_0)\wedge a)}{W^{(0)}(a)}\, \te{-(v_0-x\vee(u_0+a))\lambda(a,0)}, \qquad
x\in\mbb R.
\end{equation*}
We obtain  \eqref{eq:Min_Law} and \eqref{eq:ExpLaw} from this display
by setting $v_0=x\vee(u_0+a)$ and $u_0=-\infty$, respectively,
and recalling that $c(a)=1$.\exit

\begin{remarks}
\noindent (i)
If
$X$
is a Brownian motion with drift started at zero,
i.e.
$X_t=\mu t+\sigma W_t$,
where
$(W_t)_{t\geq0}$
is a standard Brownian motion and
$\mu\in\R, \sigma^2>0$,
then the scale function takes the form
\begin{eqnarray}
\label{eq:BM_With_Drift_Scale_Function}
W(x)& = & \frac{1}{\mu}\left[1-\te{-\frac{2\mu}{\sigma^2}x}\right]\qquad\text{if}\quad
\mu\neq0,\\
W(x)& = & \frac{2}{\sigma^2}x\qquad\text{if}\quad\mu=0.
\label{eq:BM_Scale_Function}
\end{eqnarray}
Since
$X$
has no jumps, we have an almost sure equality
$\ovl X_{\tau_a}=a+X_{\tau_a}$.
Therefore the random variable
$a+X_{\tau_a}$
is by~\eqref{eq:ExpLaw}
exponentially distributed with parameter
$(2\mu/\sigma^2)/(\te{a2\mu/\sigma^2}-1)$
(resp.
$1/a$)
if
$\mu\neq0$
(resp.
$\mu=0$).
This fact,
first observed
by Lehoczky~\cite{Lehoczky},
is generalised by
the formula in~\eqref{eq:ExpLaw}
to the class of spectrally negative L\'evy processes.\\
\noindent (ii)
Let
$X$
be a spectrally negative
$\alpha$-stable process
with cumulant generating function
$\psi(\theta)=\left(\sigma \theta\right)^\alpha$
for
$\theta>0,$
where
$\alpha\in(1,2],\,\sigma>0.$
Identity~\eqref{eq:defW}
and the definition of the Gamma function
$\Gamma$
imply that
the scale function of
$X$
takes the form
\begin{eqnarray}
\label{eq:Stable_Scale_Funct}
W(x)&=&\frac{x^{\alpha-1}}{\sigma^\alpha\Gamma(\alpha)},\qquad x\geq0.
\end{eqnarray}
Corollary~\ref{cor:Thm1}
implies that,
under $\P_x$, the random variables
$a^{-1}\le(x-\unl X_{\tau_{a-}}\ri)$
and
$a^{-1}\le(\ovl X_{\tau_a} - x\ri)$
follow
Beta$(\alpha-1,1)$
and
$\mrm{Exp}\le(\frac{\a-1}{a}\ri)$ distributions
with respective probability density functions
\begin{eqnarray*}
\P_x(x-\unl X_{\tau_{a-}}\in\td z) & = & I_{[0,a]}(z)
\frac{\alpha-1}{a} \left(\frac{z}{a}\right)^{\alpha-2}\td z,\\
\P_x(\ovl X_{\tau_a} - x \in\td z) & = & I_{(0,\infty)}(z)\frac{\alpha-1}{a}\,
\te{-z(\alpha-1)/a}\td z.
\end{eqnarray*}
Note that these distributions
do not depend on the value of
$\sigma$.
In the case of Brownian motion started at $x$
(i.e. for
$\alpha=2$), the
random variables $a^{-1}\le(x-\unl X_{\tau_a}\ri)$ and
$a^{-1}\le(\ovl X_{\tau_a} - x\ri)$
follow
$U(0,1)$
and
$\text{Exp}(1)$
distributions
respectively.
\end{remarks}

\subsection{The quadruple law}
\label{subsec:triple}
In this section we
characterise the law of the quadruple
$(\WH \t_a, \unl G_{\WH\tau_a}, \unl X_{\WH \t_a}, \ovl X_{\WH \t_a})$,
where
$\unl G_t$,
defined in~\eqref{eq:Last_time_Def},
is the last time the process
$X$
was at the infimum prior to time $t$.
Before stating the result,
we recall that $q\mapsto W^{(q)}(x)$ has a holomorphic extension
to $\C$ for every $x\geq0$ and that $(x,q)\mapsto Z^{(q)}(x)$ has
a continuous extension to $[0,\infty)\times\R$, which is
holomorphic in $q$ for every $x\in\R$ (see~\cite{P},~Lemma~2).
Furthermore for any $u\geq0$ let $\P^u$ be the exponentially
tilted probability measure defined via the Esscher transform (see
\cite[Ch. XIII]{APQ}). Then the scale functions $x\mapsto
W^{(q-\psi(u))}_u(x)$ and $x\mapsto Z^{(q-\psi(u))}_u(x)$ satisfy
the following identities for all $q\in\C, x\in\R$:
\begin{eqnarray}
\label{eq:scala_Functions_W}
W^{(q-\psi(u))}_u(x)   & = & \te{-ux} W^{(q)}(x),
\\
Z^{(q-\psi(u))}_u(x) & = & 1 + (q-\psi(u))\int_0^x\te{-uz}W^{(q)}(z)\,\td z.
\label{eq:scala_Functions_Z}
\end{eqnarray}
The identity in~\eqref{eq:scala_Functions_W}, for $q>\psi(u)$,
follows by taking Laplace transforms on both sides and
applying~\eqref{eq:defW}, and hence by analyticity for all
$q\in\C$. The identity in~\eqref{eq:scala_Functions_Z} follows
from~\eqref{eq:scala_Functions_W} and the definition
in~\eqref{eq:defZ}.

We now state
the main result of this section.

\begin{theorem}
\label{thm:Triple} Let $q, r, u, v\geq 0$. Denote $p=q+r-\psi(u)$.
Then the following identity holds:
\begin{eqnarray}
\nn\lefteqn{ \E_x\le[\te{-q\WH\tau_a - r\unl G_{\WH\tau_a}+ u\unl
X_{\WH \t_a}} I_{\left\{\ovl X_{\WH \t_a}< v\right\}}\ri]}\\
&=& \frac{W^{(q+r)}(a)}{W^{(q)}(a)}
\le[\te{-u(a-x)}\frac{Z_u^{(p)}(a+x-v)}{Z_u^{(p)}(a)} -
\te{-u(a-v)}\frac{W^{(q+r)}(a+x-v)}{W^{(q+r)}(a)}\ri].
\label{eq:Triple_Law_Equation}
\end{eqnarray}
\end{theorem}
\begin{figure}[t]
\hspace{-0.6cm}
\begin{picture}(0,0)%
\includegraphics{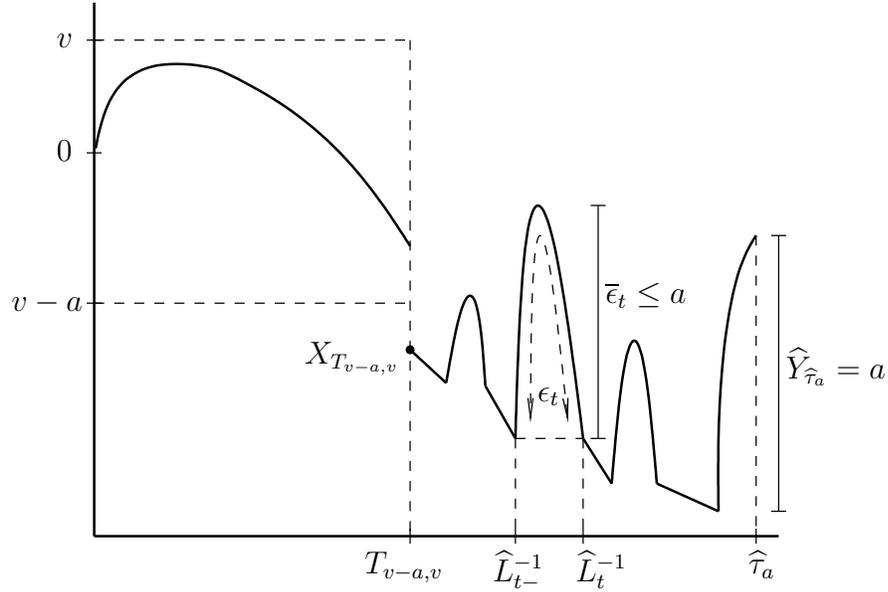}%
\end{picture}%
\setlength{\unitlength}{4144sp}%
\begingroup\makeatletter\ifx\SetFigFont\undefined%
\gdef\SetFigFont#1#2#3#4#5{%
  \reset@font\fontsize{#1}{#2pt}%
  \fontfamily{#3}\fontseries{#4}\fontshape{#5}%
  \selectfont}%
\fi\endgroup%
\begin{picture}(4665,3556)(616,-3145)
\put(2746,-3031){\makebox(0,0)[lb]{\smash{{\SetFigFont{12}{14.4}{\rmdefault}{\mddefault}{\updefault}{\color[rgb]{0,0,0}$T_{v-a,v}$}%
}}}}
\put(5041,-3031){\makebox(0,0)[lb]{\smash{{\SetFigFont{12}{14.4}{\rmdefault}{\mddefault}{\updefault}{\color[rgb]{0,0,0}$\WH \tau_a$}%
}}}}
\put(631,-1456){\makebox(0,0)[lb]{\smash{{\SetFigFont{12}{14.4}{\rmdefault}{\mddefault}{\updefault}{\color[rgb]{0,0,0}$v-a$}%
}}}}
\put(901,-556){\makebox(0,0)[lb]{\smash{{\SetFigFont{12}{14.4}{\rmdefault}{\mddefault}{\updefault}{\color[rgb]{0,0,0}$0$}%
}}}}
\put(901,119){\makebox(0,0)[lb]{\smash{{\SetFigFont{12}{14.4}{\rmdefault}{\mddefault}{\updefault}{\color[rgb]{0,0,0}$v$}%
}}}}
\put(4186,-1411){\makebox(0,0)[lb]{\smash{{\SetFigFont{12}{14.4}{\rmdefault}{\mddefault}{\updefault}{\color[rgb]{0,0,0}$\overline \epsilon_t\leq a$}%
}}}}
\put(5266,-1861){\makebox(0,0)[lb]{\smash{{\SetFigFont{12}{14.4}{\rmdefault}{\mddefault}{\updefault}{\color[rgb]{0,0,0}$\WH Y_{\WH \tau_a}=a$}%
}}}}
\put(2386,-1771){\makebox(0,0)[lb]{\smash{{\SetFigFont{12}{14.4}{\rmdefault}{\mddefault}{\updefault}{\color[rgb]{0,0,0}$X_{T_{v-a,v}}$}%
}}}}
\put(3781,-1996){\makebox(0,0)[lb]{\smash{{\SetFigFont{12}{14.4}{\rmdefault}{\mddefault}{\updefault}{\color[rgb]{0,0,0}$\epsilon_t$}%
}}}}
\put(3511,-3076){\makebox(0,0)[lb]{\smash{{\SetFigFont{12}{14.4}{\rmdefault}{\mddefault}{\updefault}{\color[rgb]{0,0,0}$\WH L^{-1}_{t-}$}%
}}}}
\put(4006,-3076){\makebox(0,0)[lb]{\smash{{\SetFigFont{12}{14.4}{\rmdefault}{\mddefault}{\updefault}{\color[rgb]{0,0,0}$\WH L^{-1}_{t}$}%
}}}}
\end{picture}
\caption{\footnotesize{A schematic representation of a path
which satisfies
$T_{v-a,v}=T_{v-a}^-$
and is therefore in the event
$\{\ovl X_{\WH \tau_a}< v\}$.
At
the moment of local time
$t$
the path in the figure has an excursion
$\epsilon_t$
away from the infimum.
The process
$\WH L^{-1}$
is the inverse local time at zero for the Markov
process
$\WH Y$.
}}
\label{fig:Min_Fig}
\end{figure}
\begin{remarks}
\noindent (i) If
$a<v$,
then
the definition of
$\WH \tau_a$
implies
$\P\left[\ovl X_{\WH \t_a}<v\right]=1$.
Therefore
Theorem~\ref{thm:Triple},
together with equations~\eqref{eq:scala_Functions_W}
and~\eqref{eq:scala_Functions_Z}
(recall that
$W^{(q)}(y)=0$
for all
$q\geq0$
and
$y<0$),
implies the identities
\begin{eqnarray}
\label{eq:Nice_Formula}
\E\left[ \te{-q\WH \t_a+ u\unl X_{\WH \t_a}} \right] & = & \te{-au}/Z_u^{(q-\psi(u))}(a),\qquad
q,u\geq0,\\
\E[\te{-q\WH\tau_a - r\unl G_{\WH\tau_a}}] & = &
\frac{W^{(q+r)}(a)}{W^{(q)}(a)Z^{(q+r)}(a)}, \qquad q,r\geq0.
\label{eq:tg}
\end{eqnarray}
The special case of  formulae~\eqref{eq:Nice_Formula}
and~\eqref{eq:tg}
for
$u=0$
and
$r=0$
respectively
(i.e. the Laplace transform of
$\WH\tau_a$)
is well-known  (see~\cite{P},~Proposition~2).
Note that in~\eqref{eq:tg}
the random variable under the expectation does not
depend on the starting point of the process
$X$.\\
\noindent (ii) Since
$X$
is spectrally negative,
we have
$\WH Y_{\WH \t_a} = a$
almost surely.
Hence formula~\eqref{eq:Triple_Law_Equation}
also yields the joint law of the
quadruple
$(\WH \t_a,\unl G_{\WH \tau_a} X_{\WH \t_a}, \ovl X_{\WH \t_a})$.\\
\noindent (iii) Under the law $\P_x$, the cdf of the random
variable $\ovl X_{\WH \t_a}-x$ is given by
\begin{eqnarray}
\label{eq:max_cdf} \P_x(\ovl X_{\WH \t_a}-x\leq z) =  1 -
\frac{W((a-z)\wedge a)}{W(a)}, \qquad z\in\R.
\end{eqnarray}
Note that the distribution of
$\ovl X_{\WH \t_a}-x$
under
$\P_x$ is absolutely continuous with respect
to the Lebesgue measure
if and only if
$W(0)=0$. If $W(0)>0$,
this distribution has an atom only at $a$,
which is of size $W(0)/W(a)$.
Furthmore, the Laplace transform of
$x-\unl X_{\WH \t_a}$
is given in terms of the scale function
$W$
of
$X$
by the formula
\begin{eqnarray}
\label{eq:Nice_Explicit_Formula}
\E_x\left[ \te{-u\left(x-\unl X_{\WH \t_a}\right)} \right] =
\te{-ua}/\left[1-\psi(u)\int_0^a\te{-uz}W(z)\,\td z\right].
\end{eqnarray}

\noindent (iv) In the special case when $X$ is a spectrally
negative $\alpha$-stable process with cumulant generating function
$\psi(\theta)=\left(\sigma \theta\right)^\alpha$,
where
$\alpha\in(1,2],\,\sigma>0,$
the scale function takes the form
$W(x)=x^{\alpha-1}/(\sigma^\alpha\Gamma(\alpha))$,
for $x\geq0$,
by~\eqref{eq:Stable_Scale_Funct}.
Formula in~\eqref{eq:max_cdf}
implies that
under $\P_x$,
for any
$x\in\R$,
the random variable
$a^{-1}\le(\ovl X_{\WH \tau_a} - x\ri)$
follows a
Beta($1,\alpha-1$) distribution
with probability distribution
\begin{eqnarray*}
\P_x(\ovl X_{\WH \tau_a} - x\in\td z) & = & I_{[0,a]}(z)\,
\frac{\alpha-1}{a} \left(1-\frac{z}{a}\right)^{\alpha-2}\td z.
\end{eqnarray*}
The Laplace transform of the random variable $a^{-1}\le(x-\unl X_{\WH \t_a}\ri)$
is, by~\eqref{eq:Nice_Explicit_Formula}, equal to
\begin{eqnarray}
\label{eq:Nice_Explicit_Formula_Stable}
\E_x\left[ \te{-ua^{-1}\left(x-\unl X_{\WH \t_a}\right)} \right] =
\te{-u}
\frac{\Gamma(\alpha)}{\Gamma(\alpha,u)},\qquad\text{where}\quad \Gamma(z,y)=\int_y^\infty \te{-s}s^{z-1}\td s
\end{eqnarray}
is the incomplete gamma function.
Note that the laws of
$a^{-1}\le(\ovl X_{\WH \tau_a} - x\ri)$
and
$a^{-1}\le(x-\unl X_{\WH \t_a}\ri)$
do not depend on the value of
$\sigma$.
In particular, in the case of Brownian motion started at $x$
(i.e.
$\alpha=2)$, the
random variables
$a^{-1}\le(\ovl X_{\WH \tau_a} - x\ri)$
and
$a^{-1}\le(x-\unl X_{\WH \t_a}\ri)$
follow $U(0,1)$
and
$\text{Exp}(1)$
distributions respectively.

\end{remarks}

\section{Application: drawdowns and rallies in exponential L\'{e}vy models}\label{sec:DD}
Let the stochastic process $S=(S_t)_{t\ge 0}$ model the price
of a stock or a foreign exchange rate. The (absolute) drawdown and rally processes of $S$ are defined by the difference $\ovl S-S$ of its running supremum $\ovl S_t = \sup_{0\leq u\leq
t}\le\{S_u\ri\}$ and the current value, and the difference $S-\unl S$ of the current
value and the running infimum $\unl S_t =
\inf_{0\leq u\leq t}\le\{S_u\ri\}$.  Their relative counterparts, the
 \textit{relative drawdown process} and \textit{relative drawup} or \textit{relative rally
processes} of $S$ are given by $\ovl S/S$ and $S/\unl S$.
Let $D_\alpha$ be the first
time the price process $S$ drops below its running supremum by at
least $(100\alpha)\%$ with $\alpha\in(0,1)$, and denote by
$U_\beta$ be the first time the stock price $S$ rallies above its
running infimum by at least $(100\beta)\%$ with $\beta>0$. In some
trading strategies buy- and sell-signals for $S$ are generated on
the basis of its relative drawdown or drawup processes. For example,
a commonly used strategy is to buy at
the epoch $U_\beta$, and to sell at the epoch $D_\alpha$
(see e.g. Levich \& Thomas~\cite{LT} for an analysis of such a strategy
in an foreign exchange setting). In this section we will consider the following four criteria
that provide descriptions of different aspects of the risk associated to investing in $S$:
\begin{itemize}
\item[(i)] The probability that $S$  attains a new (all-time) minimum at
the moment $D_\a$ of a first relative drawdown of size
$(100\a)\%$, on the event that $D_\a$ is before $T$, where $T>0$ is a given
time-horizon.
\item[(ii)] The
probability that  $S$ attains a new maximum at the
moment $U_\b$ of a first relative drawup of size $(100\beta)\%$,
on the event that $U_\b$
is before $T$.
\item[(iii)] The expected (absolute) drawdown of $S$  at the epoch $D_\a$,
on the event that $D_\alpha$ is before $T$.
\item[(iv)] The probability that a relative drawdown of size $(100\a)\%$ occurs before a
relative drawup of size $(100\beta)\%$, on the event that $D_\a$ is before
$T$.
\end{itemize}

Assume henceforth that the price process is modelled as $S=(S_t)_{t\ge 0}$ where
$$S_t=S_0\exp(X_t),\quad t\geq0,$$
with $S_0>0$ and
$X=(X_t)_{t\geq0}$
a completely asymmetric L\'evy process.
Then,
using the notation of Section~\ref{sec:exc},
we have
$\ovl S = \te{\ovl X}$,
$\unl S = \te{\unl X}$,
and therefore
$\ovl S/S=\te{Y}$,
$S/\unl S = \te{\WH Y}$,
and find
\begin{eqnarray}
\label{eq:Def_D}
D_\alpha & = &\inf\left\{t\geq0\>:\>
 \ovl S_t/S_t\geq 1/(1-\alpha)\right\}=\tau_a,\qquad\text{where}\quad
a =-\log(1-\alpha),\\
U_\beta & = & \inf\left\{t\geq0\>:\>
S_t/\unl S_t\geq 1+\beta \right\}=\WH \tau_b,\qquad\text{where}\quad
b =\log(1+\beta).
\label{eq:Def_U}
\end{eqnarray}
Furthermore, the quantities in (i)--(iv) can be expressed in the notation of Section~\ref{sec:exc}
as follows:
\begin{eqnarray*}
\P[S_{D_\alpha} = \unl S_{D_\alpha}, D_\alpha < T] &=& \P[\tau_a = \unl G_{\tau_a}, \tau_a < T],\\
\P[S_{U_\beta} = \ovl S_{U_\beta}, U_\beta < T] &=& \P[\WH\tau_b = \ovl G_{\WH\tau_b}, \WH\tau_b < T],\\
\E[(\ovl S_{D_{\alpha}}-
S_{D_\alpha})I_{\{D_\a < T\}}] &=& S_0\E[(\te{\ovl X_{\tau_a}} - \te{X_{\tau_a}})
I_{\{\tau_a< T\}}],\\
\P[D_\a < U_\beta, D_\alpha < T] &=& \P[\tau_a < \WH\tau_b, \tau_a
< T].
\end{eqnarray*}
In the following sections these quantities are explicitly expressed in terms of scale functions, employing the sextuple and quadruple
laws that were derived in Section \ref{sec:exc}.

\subsection{The probability of attaining a new minimum at the first moment of a drawdown}
In the next result the distributions of the two random vectors $(\tau_a, \WH Y_{\tau_a})$
and $(\WH\tau_a, Y_{\WH\tau_a})$ are explicitly identified in terms of scale functions.
Note that the process $X$ is equal to its running minimum at the first moment of
drawdown of size $a$ if and only if $\WH Y_{\tau_a}$ is equal to zero.
Similarly, $X$ is equal to its running minimum at the first moment
of a drawup of size $b$ if and only if $Y_{\WH\tau_a}$ is equal to zero.


\begin{cor}
\label{thm:JointLawT_tau} Let $a>0$. (i)
For $b, q,\theta\geq0$ the following identities hold true:
\begin{align}
\label{eq:Joint_Law_T_tau} \E\left[\te{-q\tau_a}I_{\{\WH
Y_{\tau_a} > b\}}\right] & =
\te{-b\lambda(a,q)} \left[
1 - \lambda(a,q) \int_0^a
\frac{W^{(q)}(y)}{W^{(q)}(a)}\,\td y \right] M_{q,a},\\
\E\left[ \te{-q\WH \t_a-\theta Y_{\WH \t_a}} \right]  &= 1-\left(
\frac{(q-\psi(\theta))\te{-a\theta}}{Z_\theta^{(q-\psi(\theta))}(a)}+\frac{\theta}{W^{(q)}(a)}
\right) \int_0^aW^{(q)}(z)\,\td z, \label{eq:Laplace_Transf_joint}
\end{align}
where $\lambda({a,q}) = W^{(q)\prime}_+(a)/W^{(q)}(a)$ and
$M_{q,a}:= \E[\te{-\lambda(a,q)(Y_{\tau_a}-a)}]$ is given by
\begin{equation}\label{eq:M}
M_{q,a}  = Z_\lambda^{(-\psi(\lambda))}(a) -
W_\lambda^{(-\psi(\lambda))}(a)\frac{-\psi(\lambda)W_\lambda^{(-\psi(\lambda))}(a)+\lambda Z_\lambda^{(-\psi(\lambda))}(a)}
{W_{\lambda+}^{(-\psi(\lambda))\prime}(a)+\lambda W_\lambda^{(-\psi(\lambda))}(a)}\quad \text{with}\>\> \lambda=\lambda(a,q).
\end{equation}

(ii) The Laplace transforms of
$t\mapsto \P(\tau_a < t, \tau_a = \unl G_{\tau_a})$
and $t\mapsto \P(\WH\tau_a < t, \WH\tau_a = \ovl G_{\WH\tau_a})$
are given by
\begin{eqnarray}
\label{eq:Y0q}
\int_0^\infty \te{-qt} \P[\tau_a < t, \tau_a = \unl G_{\tau_a}]\td t &= &
\frac{1}{q}\left\{ 1 - \lambda(a,q) \int_0^a
\frac{W^{(q)}(y)}{W^{(q)}(a)}\,\td y\right\}  M_{q,a},\\
\int_0^\infty \te{-qt} \P[\WH\tau_a < t, \WH\tau_a = \ovl G_{\WH\tau_a}]\td t &=&
\frac{1}{q}\left\{ 1 - \lambda(a,q)\int_0^a \frac{W^{(q)}(y)}{W^{(q)}(a)}\td y\right\}.
\label{eq:Yh0q}
\end{eqnarray}
\end{cor}

\begin{remarks}
\noindent (i) From the formulas in~\eqref{eq:Y0q} and \eqref{eq:Yh0q}
it follows that the non-negative random variables $\WH Y_{\tau_a}$ and
$Y_{\WH\tau_a}$
have atoms at zero of the sizes
\begin{eqnarray}
\label{eq:atom_at_zero_wh_Y} \P\left[\WH Y_{\tau_a}=0\right] &=& 1 -
\le(1 - \frac{W_+'(a)}{W(a)^2} \int_0^a W(y)\td y\ri)
M_{0,a},\\
\label{eq:YH0}
\P\left[Y_{\WH \t_a}=0\right]  &=&  1 -
\frac{W_+'(a)}{W(a)^2} \int_0^aW(y)\,\td y.
\end{eqnarray}
Furthermore, the formula in~\eqref{eq:Joint_Law_T_tau}, implies that
conditional on the event $\left\{\WH
Y_{\tau_a}>0\right\}$,
$\WH Y_{\tau_a}$
follows an exponential distribution with mean
$W(a)/W_+'(a)$.

\noindent (ii) In the case $X$ is Brownian motion with non-zero
drift $\mu$ and Gaussian coefficient $\sigma^2$, the scale
function is given by~\eqref{eq:BM_With_Drift_Scale_Function}.
Hence the atoms of $\WH Y_{\tau_a}$ and $Y_{\WH\tau_a}$ at zero
are by~\eqref{eq:atom_at_zero_wh_Y} and \eqref{eq:YH0} of the
sizes
$$
\P\left[\WH Y_{\tau_a}=0\right]  =
\frac{\te{-a2\mu/\sigma^2}-1+a2\mu/\sigma^2}{(\te{-a\mu/\sigma^2}-\te{a\mu/\sigma^2})^2}.
\qquad \P\left[Y_{\WH \tau_a}=0\right]  =
\frac{\te{a2\mu/\sigma^2}-1-a2\mu/\sigma^2}{(\te{-a\mu/\sigma^2}-\te{a\mu/\sigma^2})^2}.
$$
Furthermore, $\WH Y_{\tau_a}$ and $Y_{\WH\tau_a}$ conditioned to
be stirctly positive follow exponential distributions:
$$
Y_{\WH\tau_a}\left|Y_{\WH\tau_a}\right.>0 \sim
\text{Exp}\left(\frac{2\mu/\sigma^2}{\te{2a\mu/\sigma^2}-1}\right),
\qquad \WH Y_{\tau_a}\left|\WH Y_{\tau_a}>0\right.\sim
\text{Exp}\left(\frac{2\mu/\sigma^2}{1-\te{-2a\mu/\sigma^2}}\right).
$$
If $\mu=0$, formula~\eqref{eq:atom_at_zero_wh_Y} implies $
\P\left[\WH Y_{\tau_a}=0\right]  = 1/2$. In this case
formula~\eqref{eq:Joint_Law_T_tau} implies that $\P\left[\WH
Y_{\tau_a}\in \td b|\WH Y_{\tau_a}>0\right]  = (\te{-b/a}/a)\td b$
for $b>0$. These results coincide with the findings of~\cite[Prop.
2.2,\ Cor. 2.6]{PosVec_1}.

(i) For later reference we note that the
formula in~\eqref{eq:Laplace_Transf_joint} can
equivalently be expressed as
\begin{eqnarray}
\label{eq:Equiv_Form_For_LAplace_Pair} \E\left[ \te{-q\WH
\t_a-\theta Y_{\WH \t_a}} \right] & = & 1-
\frac{\int_0^\infty\te{-\theta z}W_+^{(q)\prime}(a+z)\,\td
z}{\int_0^\infty\te{-\theta z}W^{(q)}(a+z)\,\td z}\,
\int_0^a\frac{W^{(q)}(z)}{W^{(q)}(a)} \,\td z
\end{eqnarray}
for any $q\geq0$ and $\theta>\Phi(q)$.
\end{remarks}

\smallskip
\noindent {\it Proof }{\it of Corollary \ref{thm:JointLawT_tau}}
The identity 
$ \WH Y_{\tau_a} = \ovl X_{\tau_a} -  Y_{\tau_a} - \unl X_{\tau_a} $
implies that, for any $h>0$, on the event
$$\left\{\WH Y_{\tau_a}\leq b,Y_{\tau_a}-a\in\td h\right\}=\left\{
\ovl X_{\tau_a}\leq a+b+h +\unl X_{\tau_a}, Y_{\tau_a}-a\in\td
h\right\}$$ the inequality $\ovl X_{\tau_a}\leq a+b+h $ holds
$\P_0$-a.s.
%
The formula in~\eqref{eq:QuitLaq_overshoot} of
Theorem~\ref{thm:lawreflect1} and the definition of
$F_{q,q,a}$ given in~\eqref{Def:F} imply
\begin{eqnarray}\label{e1}
\E\left[\te{-q\tau_a}I_{\{\WH Y_{\tau_a}\leq b, Y_{\tau_a}-a\in\td h\}}\right]  = I_1(\td h)
+ I_2(\td h),
\end{eqnarray}
where $I_1$ and $I_2$ are measures on $(\mbb R_+,\mc B(\mbb R_+))$ given by
\begin{eqnarray*}
I_1(\td h) &:=& \E\left[\te{-q\tau_a}I_{\{\WH Y_{\tau_a}\leq b, Y_{\tau_a}-a\in\td h,\ovl X_{\tau_a}
\leq b+h\}}\right] = \nu_{q,a}^0(\td h)- \te{-b\lambda(a,q)}
\nu_{q,a}^1(\td h),\\
I_2(\td h) &:=& \E\left[\te{-q\tau_a}I_{\{\WH Y_{\tau_a}\leq b, Y_{\tau_a}-a\in\td h,\ovl X_{\tau_a}\in
(b+h, a + b+h]\}}\right]  = \lambda(a,q)\te{-b\lambda(a,q)}\nu_{q,a}^1(\td h)
\int_{0}^a\frac{W^{(q)}(y)}{W^{(q)}(a)}\,\td y.
\end{eqnarray*}
In these expressions 
$\lambda(a,q)$ is given in \eqref{eq:laq} and
$\nu_{q,a}^j$, $j\in\{0,1\}$, are measures supported on
$\mbb R_+$ and defined by the following formula
\begin{eqnarray}
\label{eq:Nu_a_Def} \nu_{q,a}^j\left(\td h\right) =
\te{-hj\lambda(a,q)}\int_0^a
R^{(q)}_{a}(\td y) \Lambda\left(y-a-\td h\right),
\end{eqnarray}
where $R^{(q)}_a$ is the resolvent measure  given in~\eqref{eq:res} and
$\Lambda$ is the L\'evy measure of $X$.
An analogous argument, based on the
formula in~\eqref{eq:QuitLaq_Creeping} of
Theorem~\ref{thm:lawreflect1}, yields the formula
\begin{equation}\label{e2}
\E\left[\te{-q\tau_a}I_{\{\WH Y_{\tau_a}\leq b,
Y_{\tau_a}=a\}}\right] = \le[1- \te{-b\lambda(a,q)} +  \lambda(a,q)\te{-b\lambda(a,q)}\int_{0}^a\frac{W^{(q)}(y)}{W^{(q)}(a)}\,\td y\ri] \Delta^{(q)}(a),
\end{equation}
where $\Delta^{(q)}(a)$ is given in \eqref{eq:Daq}.
Observe that, on account of the
compensation formula applied to $(\Delta X_t)_{t\ge 0}$, the
following relations hold true:
\begin{eqnarray}
\label{eq:mq} \E[\te{-q\tau_a}] & = & \nu_{q,a}^0(0,\infty)+
\Delta^{(q)}(a),\\
\E[\te{-\lambda(Y_{\tau_a}-a)}]
& =  & \nu_{q,a}^1(0,\infty) + \Delta^{(q)}(a), \q \lambda=\lambda(a,q),\label{eq:Mq}
\end{eqnarray}
where the expressions for $\E[\te{-q\tau_a}]$ and $\E[\te{-\lambda(Y_{\tau_a}-a)}]$ in terms of scale functions
are given in \eqref{eq:taua} and \eqref{eq:M} (see~\cite[Thm. 1]{AKP}).
By integrating \eqref{e1}
over $h\in(0,\infty)$, adding \eqref{e2} and using \eqref{eq:mq}, \eqref{eq:Mq} and
\eqref{eq:taua} we obtain the formula in \eqref{eq:Joint_Law_T_tau}.

It follows from Theorem \ref{thm:Triple} that the identity
\begin{eqnarray}
\E\le[\te{-q\WH\tau_a + u X_{\WH \t_a}} I_{\left\{\ovl X_{\WH \t_a}< v\right\}}\ri]
= \frac{Z_u^{(q-\psi(u))}(a-v)}{Z_u^{(q-\psi(u))}(a)} -
\te{uv}\frac{W^{(q)}(a-v)}{W^{(q)}(a)}.
\label{eq:Laplace_Joint_Law}
\end{eqnarray}
holds for $q,u,v\geq0$. The Laplace transform in $v$ of the
identity in~\eqref{eq:Laplace_Joint_Law} evaluated at $u=\theta$
implies the equality in~\eqref{eq:Laplace_Transf_joint}. 
Since
$\E\left[\te{-q\WH\tau_a}I_{\{Y_{\WH \t_a}=0\}}\right] = \lim_{\theta\to\infty}\E\left[
\te{-q\WH\tau_a -\theta Y_{\WH \t_a}} \right]$, we obtain the formula in \eqref{eq:Yh0q}
by noting that $\{Y_{\WH \t_a}=0\}= \{\WH\tau_a = \ovl G_{\WH\tau_a}\}$ and applying the identity
$$
\lim_{\theta\to\infty}\frac{\theta\int_0^\infty\te{-\theta z}W_+^{(q)\prime}(a+z)\,\td z}%
{\theta\int_0^\infty\te{-\theta z}W^{(q)}(a+z)\,\td z} =
\frac{W^{(q)\prime}_+(a)}{W^{(q)}(a)},
$$
which follows from well-known properties of the scale functions.
\exit

\subsection{The expected drawdown $\ovl S-S$ at $D_\a$.}

\begin{cor}
The Laplace transform of $T\mapsto \E[(\ovl S_{D_\alpha}-S_{D_\alpha})\mbf 1_{\{D_\alpha < T\}}]$
is given by
\begin{eqnarray}\nn
\lefteqn{\int_0^\infty \te{-qT} \E[(\ovl S_{D_\alpha} - S_{D_\alpha})\mbf 1_{\{D_\alpha < T\}}]\td T} \\
&=&
\frac{S_0\lambda}{\lambda - 1}\left[\int_0^a h_q(x)W^{(q)}(x)\td x\right]
- \frac{S_0}{\lambda - 1}
\left[ \int_0^a h_q(x)W^{(q)}(\td x) - \frac{1}{q}\right], \qquad q>\psi(1),
\label{eq:oversh}
\end{eqnarray}
where $a=-\log(1-\alpha)$,
$\lambda = \lambda(a,q) = W^{(q)\prime}_+(a)/W^{(q)}(a)$ and $h_q(x) =1 -  (1-q^{-1}\psi(1))\te{-x} $.
\end{cor}
\proof Recall that $a =-\log(1-\alpha).$
In the notation of Section \ref{sec:exc}, the Laplace transform on the left-hand side of
equation~\eqref{eq:oversh} is equal to
\begin{equation*}
\frac{S_0}{q} \E[\te{-q\tau_a}\{\te{\ovl X_{\tau_a}} - \te{X_{\tau_a}}\}]
= \frac{S_0}{q} \E[\te{-q\ovl G_{\tau_a} + \ovl X_{\tau_a}}\{\te{-q(\tau_a - \ovl G_{\tau_a})}
(1 - \te{-Y_{\tau_a}})\}].
\end{equation*}
Direct integration of the expressions in Theorem~\ref{thm:lawreflect1} with $u=-\infty$ and $x=0$ shows that,  for $q>\psi(1)$,
\begin{eqnarray}\label{eq:Lt0}
\lefteqn{\E[\te{-q\tau_a}\{\te{\ovl X_{\tau_a}} - \te{X_{\tau_a}}\}]}\\
&=&  \E[\te{-q\ovl G_{\tau_a} + \ovl X_{\tau_a}}]\{\E[\te{-q(\tau_a - \ovl G_{\tau_a})}]
- \E[\te{-q(\tau_a - \ovl G_{\tau_a})-Y_{\tau_a}}]\},\nn
\end{eqnarray}
with
\begin{eqnarray}\label{eq:Lt1}
\E[\te{-q(\tau_a - \ovl G_{\tau_a})-(Y_{\tau_a}-a)}] &=&
\frac{\lambda(a,q)}{\lambda(a,0)}\nonumber
\\ &\phantom{=}& \cdot \le[Z_1^{(q-\psi(1))}(a) - W_1^{(q-\psi(1))}(a)
\frac{(q-\psi(1))W_1^{(q-\psi(1))}(a) +
Z_1^{(q-\psi(1))}(a)}{W_{1+}^{(q-\psi(1))\prime}(a)+W_1^{(q-\psi(1))}(a)}\ri],\\
\label{eq:Lt2}
\E[\te{-q(\tau_a - \ovl G_{\tau_a})}] &=&
\frac{\lambda(a,q)}{\lambda(a,0)}
\le[Z^{(q)}(a)
-q \lambda(a,q)^{-1} W^{(q)}(a)\ri],\\
\E[\te{-q\ovl G_{\tau_a} + \ovl X_{\tau_a}}] &=&  \frac{\lambda(a,0)}{\lambda(a,q)}
\cdot \frac{\lambda(a,q)}{\lambda(a,q) - 1},
\end{eqnarray}
where the functions $Z_1^{(q-\psi(1))}(a)$ and $W_1^{(q-\psi(1))}$ are defined in 
equations~\eqref{eq:scala_Functions_Z} and~\eqref{eq:scala_Functions_W}. Inserting the expressions in equations~\eqref{eq:Lt1} and \eqref{eq:Lt2}
into equation~\eqref{eq:Lt0} yields the expression in equation~\eqref{eq:oversh}. \exit

\subsection{Probability of a large drawdown preceding or following a small rally}
\label{sec:Drawdown_and_Drawup} In this section we identify explicit expressions for the probabilities of the events $\{D_\alpha
<U_\beta\}$ and $\{D_\alpha >U_\beta\}$ in the case when the size
$1/(1-\alpha)$ of the relative drawdown is larger or equal to the
size $1+\beta$ of the relative drawup. Note that the case
$1/(1-\alpha)\geq1+\beta$ is the case of most interest in practice,
since the event $\{D_\alpha <U_\beta\}$ corresponds to the risk
held by an investor of losing more in an investment than
they would gain if the security rallied. We refer to Hadjiliadis
\& Vecer~\cite{PosVec_1} for further background on the use of this
probability as risk-measure.

The results in this section are a consequence of
Corollary~\ref{thm:JointLawT_tau} and the following key identity:

\begin{lem}
\label{lem:TrnasformEvent}
Pick any numbers
$a\geq b>0$. Then the following equality holds almost surely:
\begin{eqnarray}
\label{eq:drawdown_Frist}
\left\{\tau_a<\WH \tau_b\right\}
=
\left\{a-b<Y_{\WH \tau_b}\right\}.
\end{eqnarray}
\end{lem}

\begin{remark}
In the case
$b>a$
and the L\'evy measure of
$X$
is non-zero,
the equality analogous to~\eqref{eq:drawdown_Frist}
does not hold.
The event
$\left\{\tau_a<\WH \tau_b\right\}$ 
is different from $\{b-a>\WH Y_{\tau_a}\}$ as it depends in an
essential way on the supremum $\ovl{\WH Y}_{\tau_a}$ of the
reflected process $\WH Y$ at $\tau_a$. The identification of the
law of the random variable $\ovl{\WH Y}_{\tau_a}$ is beyond the
scope of the current paper.
\end{remark}

\smallskip
\noindent {\it Proof }{\it of Lemma \ref{lem:TrnasformEvent}.} Let
$\ovl Y$ and $\ovl {\WH Y}$ denote the processes given by the
running suprema of the reflected processes $Y$ and $\WH Y$
respectively. In Hadjiliadis \& Vecer~\cite{PosVec_1} the
following model-free equalities were established:
\begin{eqnarray}
\label{eq:Key_Equality_Fro_LEmma} Y+\WH Y =   \ovl Y \vee \ovl
{\WH Y} = \ovl {\WH Y} +\left[0\vee \left(\ovl Y - \ovl {\WH Y}
\right)\right]
\end{eqnarray}
Since
$\ovl {\WH Y}_{\WH \tau_b}=\WH Y_{\WH \tau_b}=b$ a.s.,
we find $\left\{\tau_a<\WH \tau_b\right\}   = \{a < \ovl Y_{\WH\tau_b}\} =
\left\{a-b< \ovl Y_{\WH \tau_b}-\ovl{\WH Y}_{\WH \tau_b}\right\}$,
and the lemma follows in view of~\eqref{eq:Key_Equality_Fro_LEmma}.
\exit

By combining Lemma~\ref{lem:TrnasformEvent} with
Corollary~\ref{thm:JointLawT_tau} we obtain a semi-analytical
expression for the probability of interest:

\begin{cor}
\label{thm:CondLaw}
Let $a\geq b>0$ and $t\ge 0$.
Then we have
$\P\left[\WH \tau_b<\tau_a\wedge t\right] =  \P\left[Y_{\WH \tau_b}\leq a-b, \WH \tau_b<t\right],$
where the joint Laplace transform of
$(u,t)\mapsto\P\left[Y_{\WH \tau_b}\leq u, \WH \tau_b< t \right]$
is given by the formula
\begin{eqnarray}
\nn\lefteqn{\int_0^\infty\int_0^\infty \te{-\theta u - q t}
\P\left[Y_{\WH \tau_b}\leq u, \WH \tau_b < t \right]\td u = }\\
&& \frac{1}{q\theta}\left[1 +\left(
\frac{(\psi(\theta)-q)\te{-b\theta}}{Z_\theta^{(q-\psi(\theta))}(b)}-\frac{\theta}{W^{(q)}(b)}
\right) \int_0^bW^{(q)}(y)\,\td y\right].
\label{eq:Up_before_Down}
\end{eqnarray}
for any
$\theta,q > 0$.
Moreover for any
$a>0$
we have
\begin{eqnarray*}
\P\left[\tau_a<\WH \tau_a\right] & = &  \frac{W_+'(a)}{W(a)^2}\int_0^aW(z)\,\td z.
\end{eqnarray*}
\end{cor}

\begin{remarks}
\noindent (i)
It is a consequence of the representation of
the scale function
$Z_\theta^{(-\psi(\theta))}$
given in~\eqref{eq:scala_Functions_Z},
that the Laplace transform in~\eqref{eq:Up_before_Down}
can be expressed as
$$
\theta \int_0^\infty \te{-\theta u}
\P\left[Y_{\WH \tau_b}\leq u \right]\td u
 =
1+\left(
1/\int_0^\infty\te{-z\theta} W(b+z)\,\td z - \theta/ W(b)
\right)
\int_0^bW(y)\,\td y.
$$
\noindent (ii) The probability measure in
Corollary~\ref{thm:CondLaw} does not depend on the starting point
of the process $X$, since the reflected processes $Y, \WH Y$, and
hence the stopping times $\tau_a, \WH \tau_b$, are independent of
$X_0$. Therefore equality in~\eqref{eq:Up_before_Down} remains
valid if $\P[\cdot]$ is replaced by $\P_x[\cdot]$ for any starting
point
$x\in\R$.\\
\noindent (iii) Recall
from~\eqref{eq:BM_With_Drift_Scale_Function} that in the case $X$
is a Brownian motion with drift $\mu\in\R\setminus\{0\}$ and the
Gaussian coefficient $\sigma^2>0$, the scale function takes the
form $W(x)=\left(1-\te{-x2\mu/\sigma^2}\right)/\mu$ for $x\geq0$.
The formula in Corollary~\ref{thm:CondLaw} yields
$$
\P\left[\tau_a<\WH \tau_a\right] =
\frac{\te{-a2\mu/\sigma^2}-1+a2\mu/\sigma^2}{\left(\te{a\mu/\sigma^2}-\te{-a\mu/\sigma^2}\right)^2},
$$
which coincides with the formula in~(2.1) of~\cite{PosVec_1}. In
particular if $\mu=0$ then Corollary~\ref{thm:CondLaw} yields
$\P\left[\tau_a<\WH \tau_a\right] =1/2$, which follows of course
also directly by symmetry.
\end{remarks}

\subsection{Example: Carr \& Wu model}
\label{sec:cas}

Carr \& Wu~\cite{CarrWu} documents the persistence of
the implied volatility skew across maturities
in the exchange traded options on the S\&P 500 index.
The evidence presented in~\cite{CarrWu}
suggests that the
left tail of the risk-neutral
distribution of returns of the index remains ``fat''
as the maturity increases.
In order to capture this phenomenon
Carr \& Wu~\cite{CarrWu}
model the risk-neutral evolution of the
index level by the following
stochastic differential equation:
\begin{eqnarray}
\label{eq:MOdel_MarketCrash}
\td S_t/S_{t-}  = (r-q)\,\td t+\Sigma\, \td L_t^{\alpha,-1}, \qquad \alpha\in(1,2),\,\Sigma>0,
\end{eqnarray}
where
$r$
and
$q$
denote the continuously compounded risk free rate and dividend yield. The
driver
$L^{\alpha,-1}$
is assumed to be an
$\alpha$-stable L\'evy martingale with maximal negative skewness, i.e.
$L^{\alpha,-1}$
is
a spectrally negative
$\alpha$-stable process.
In this model the decay of the left tail of the log return
(i.e. of the quantity
$\P\left[\log(S_T/S_t)<s\right]$
for
any fixed
$0\leq t<T$)
is asymptotically equal to
$|s|^{-\alpha}$
as
$s$
tends to
$-\infty$,
while the right tail decays exponentially, so that
the model~\eqref{eq:MOdel_MarketCrash}
captures the observed phenomenon of the persistence of
the implied volatility skew across maturities.

In order to understand the risk-neutral probabilities of drawdowns and rallies in this model
for the S\&P~500 index,
note that the process
$S$
can be expressed as
$S=\exp(X)$,
where
$X$
is
a spectrally negative
$\alpha$-stable process
with drift.
The cumulant generating function
of
$X$
takes the form
\begin{equation}
\label{eq:stable_drift_cumulant}
\psi(\theta)=\mu \theta+\left(\sigma \theta\right)^\alpha\quad\text{for}\quad \theta>0,\quad \text{where}\quad
\alpha\in(1,2),
\end{equation}
for some
$\mu\in\R,\sigma>0$
that depend on
$r,q,\Sigma$
and
$\alpha$.
Note that in the limit case of
$\alpha=2$,
the process
$X$
corresponds to Brownian motion with drift
$\mu$
and variance
$2\sigma^{2}$
at time one,
and in this case the model~\eqref{eq:MOdel_MarketCrash}
reduces to the Black-Scholes model.

As can be seen from Proposition~\ref{thm:CondLaw},
the probability of a drawdown preceding a rally in
the model~\eqref{eq:MOdel_MarketCrash}
is given in terms of scale functions of the process
$X$,
which is expressed in terms of the
Mittag-Leffler function
$$
E_{\beta,\gamma}(y)=\sum_{n=0}^\infty\frac{y^n}{\Gamma(n\beta+\gamma)},\quad y\in\R.
$$
Note that $E_{\beta,\gamma}$ is an entire function if
$\beta,\gamma$ are strictly positive and that it generalises the exponential
$E_{1,1}(y)=\te{y}$ and the hyperbolic cosine
$E_{2,1}(y^2)=\cosh(y)$ functions. Furrer~\cite{HFurrer}
identified the 0-scale function $W$ of $X$ as
\begin{eqnarray}
\label{eq:stable_drift_scale_fun} W(x) = \frac{1}{\mu}\left[1-
E_{\alpha-1,1}\left(-\frac{\mu}{\sigma^\alpha}x^{\alpha-1}\right)\right],\qquad
x\geq0, \a\in(1,2].
\end{eqnarray}
As a consequence, the $q$-scale function $W^{(q)}$ of $X$, which
is related to $W$ by the well-known expression $W^{(q)}(x) =
\sum_{k\ge 0} q^k W^{\star(k+1)}(x)$ that is given in terms of the $k$-th
convolution $W^{\star k}$, $k\in\mbb N$, of $W$ with itself,
admits the following series representation:
\begin{equation}
W^{(q)}(x) = \frac{1}{\mu}\le[1 - E_{\a-1,1,1}\le(- \frac{\mu}{\s^\a}x^{\a - 1},
-\frac{qx}{\mu}\ri)\ri], \qquad x,q\ge 0, \a\in(1,2],
\end{equation}
where, for any nonnegative $\b$, $\g$ and $\d$,
 the function $E_{\b,\g,\d}: \mathbb R^2\to\mathbb R$ is the function that
is defined by
$$
E_{\b,\g,\d}(y,z) :=  \sum_{n=0}^\infty\sum_{k=0}^\infty
\frac{s_{k,n}}{\Gamma(n\b+ k\d +\g )} z^k y^n,
\q y,z\in\mbb R,
$$
where the coefficients $(s_{k,n}, n,k\in\mbb N\cup\{0\})$ are given by
\begin{eqnarray*}
s_{k,n+1} = \le(\begin{array}{c}{n}\\ k\end{array} \ri), \qquad k, n \in \mbb N\cup\{0\}, \qquad s_{0,0}=1.
\end{eqnarray*}
Note that $s_{k,n}=\frac{(n-1)!}{k! (n-1-k)!}$ if $k=0, 1, \ldots, n-1$,
and $s_{k,n}$ is equal to zero otherwise.

\begin{remarks}
\noindent (i) In the limit as
$\alpha$
approaches $2$,
the formula in~\eqref{eq:stable_drift_scale_fun}
tends to the scale function of the Brownian motion with volatility
$\sigma\sqrt{2}$,
which is given in~\eqref{eq:BM_With_Drift_Scale_Function}.\\
\noindent (ii) It follows from the formula in~\eqref{eq:stable_drift_scale_fun}
that if the modulus of the drift
$\mu$
tends to zero, then
in the limit we obtain the scale function,
of the
spectrally negative
$\alpha$-stable process,  which is given in~\eqref{eq:Stable_Scale_Funct}.\\
\noindent (iii) It is assumed in~\eqref{eq:stable_drift_scale_fun}
 that $\alpha>1$ and
hence the series $E_{\alpha-1,1}$ converges uniformly on compact
subsets. Therefore the derivative
and the integral
of the scale function
can be obtained by performing these operations under the summation:
\begin{eqnarray}
\label{eq:stable_drift_scale_fun_derivative}
W'(x) & = &\frac{\alpha-1}{\sigma^\alpha}x^{\alpha-2} E_{\alpha-1,1}'\left(-\frac{\mu}{\sigma^\alpha}x^{\alpha-1}\right),\qquad
x>0.\\
\int_0^aW(z)\,\td z & = &
\frac{a}{\mu}\left[1- E_{\alpha-1,2}\left(-\frac{\mu}{\sigma^\alpha}a^{\alpha-1}\right)\right],\qquad a>0.
\label{eq:scale_fn_int}
\end{eqnarray}
\end{remarks}

\smallskip

We now present an analytically tractable expression for the
probability of a large drawdown occurring
before a small rally.
Let
$\mathcal L^{-1}$
denote the inverse Laplace transform.
\begin{cor}
\label{cor:Crash_Dif}
Let
$a=-\log(1-x)$,
$x\in(0,1)$,
and
$b=\log(1+y)$,
$y>0$,
and
$a>b$.
Then the probability of
a drawdown of size
$(100x)\%$,
in the Carr \& Wu  model~\eqref{eq:MOdel_MarketCrash}
occurs before a rally of size
$(100y)\%$
is given by
$\P\left[\tau_a<\WH \tau_b\right] = \mathcal L^{-1}(F)(a-b)$
where
\begin{equation}
\label{eq:Laplace_Formula_Drawdawn}
F(\theta) =
\frac{1}{\theta}-\left[
\mu\te{-b\theta}/
\left(\sum_{n=1}^\infty
\left(-\frac{\mu}{\sigma^\alpha\theta^{\alpha-1}}\right)^n
\frac{\Gamma(n(\alpha-1)+1,b\theta)}{\Gamma(n(\alpha-1)+1)}
\right)
+1/W(b)
\right]
\int_0^bW(z)\,\td z,
\end{equation}
where
$\Gamma(n(\alpha-1)+1,b\theta)$
denotes the incomplete gamma function (see~\eqref{eq:Nice_Explicit_Formula_Stable})
and the scale function $W$
and its anti-derivative are explicitly given in~\eqref{eq:stable_drift_scale_fun}
and~\eqref{eq:scale_fn_int}.
\end{cor}

\begin{remark}
Standard Laplace inversion algorithms, such as the one given in
Abate and Whitt~\cite{AbateWhitt},
would typically evaluate the function in~\eqref{eq:Laplace_Formula_Drawdawn}
about twenty times in order to compute the probability
$\P\left[Y_{\WH \tau_b}\leq a-b \right]$.
Since the series in expression~\eqref{eq:Laplace_Formula_Drawdawn}
is dominated by a geometric series,  the error at each
evaluation of the Laplace transform can be controlled easily.
\end{remark}

In the case when the relative drawdown level $1/(1-x)$ is strictly
larger than the size $1+y$ of the relative drawup (i.e. rally) of
the index $S$, the probability of the market crash
$\P\left[\tau_a<\WH \tau_b\right]$, where $a=-\log(1-x)$ and
$b=\log(1+y)$, is by Proposition~\eqref{thm:CondLaw} equal to
$\P\left[Y_{\WH \tau_b}\leq a-b \right]$. By
Proposition~\ref{thm:CondLaw} the Laplace transform of the
function $u\mapsto \P\left[Y_{\WH \tau_b}\leq u \right]$ is given
by the formula in~\eqref{eq:Laplace_Formula_Drawdawn}.

The previous results yield a closed-form expression
for the probability of
a drawdown occurring before a rally in the case these are of the
same size in the log scale (cf.~\eqref{eq:Def_D} and~\eqref{eq:Def_U}):
\begin{cor}
\label{cor:Prob_a_a}
The probability of
a drawdown of size
$(100x)\%$,
$x\in(0,1)$,
in the Carr \& Wu  model~\eqref{eq:MOdel_MarketCrash}
occurs before a rally of size
$(100x/(1-x))\%$
is given by
\begin{eqnarray}
\label{eq:Formula_Crash_Symmetric}
\P\left[\tau_a<\WH \tau_a\right] = (\alpha-1)A E'_{\alpha-1,1}\left[-A\right]
\frac{1-E_{\alpha-1,2}\left[-A\right]}{\left(1-E_{\alpha-1,1}\left[-A\right]\right)^2},\quad
\text{with}\quad A = \frac{\mu}{\sigma^\alpha}a^{\alpha-1},
\end{eqnarray}
where
$a=-\log(1-x)$.
\end{cor}
Note that in the limiting case
$\alpha=2$,
formula~\eqref{eq:Formula_Crash_Symmetric}
yields the well-known expression for Brownian motion with drift $\mu$
and volatility
$\sigma\sqrt{2}$.
Formula~\eqref{eq:Formula_Crash_Symmetric}
follows from
Proposition~\ref{thm:CondLaw}
and the formulas~\eqref{eq:stable_drift_scale_fun},~\eqref{eq:stable_drift_scale_fun_derivative}
and~\eqref{eq:scale_fn_int}.

\section{Proofs}
\label{sec:proofs}

\subsection{Proof of Theorem~\ref{thm:lawreflect1}}
We now give a proof of Theorem~\ref{thm:lawreflect1}
based on It\^o's excursion theory.

Start by noting that since both sides
in~\eqref{eq:QuitLaq_overshoot}
and~\eqref{eq:QuitLaq_Creeping}
are continuous in
$q$,
it suffices to prove the formulae
in~\eqref{eq:QuitLaq_overshoot}
and~\eqref{eq:QuitLaq_Creeping}
for
$q>0$,
which will be assumed without loss of generality.
We prove the result in three steps.
Step~(1) deals with the segment of a path in $$A_o=\left\{ \unl
X_{\tau_a}\geq u, \ovl X_{\tau_a}\in \td v, Y_{\tau_a-}\in\td y,
Y_{\tau_a}-a\in\td h \right\}$$ over the time interval
$[0,T_{u,u+a}]$ by applying the strong Markov property at time
$T_{u,u+a}$. Step~(2) extracts the final jump $\Delta X_{\tau_a}$
in the expectation of~\eqref{eq:QuitLaq_overshoot} by applying the
compensation formula for the Poisson point process $(\Delta
X_t)_{t\geq0}$. Step~(3) applies the It\^o's excursion theory to
the segment of a path in $A_o$ over the interval
$(T_{u,u+a},\tau_a)$. In the case of the event $A_c$ the structure
of the proof is similar with steps~(2) and~(3) merged as there is
no overshoot at time $\tau_a$.

\textit{Step (1)} Note that
$\P_x\left[\left\{T_{u,u+a}=T^-_u\right\}\cap A_o\right] =  0$
(see Figure~\ref{fig:Main})
and that on
$A_o$
we have
$T_{u,u+a}<\ovl G_{\tau_a}$
$\P_x$-a.s.
We can therefore
apply the strong Markov property at
$T_{u,u+a}$
and identity~\eqref{eq:two-sided1}
to the expectation in~\eqref{eq:QuitLaq_overshoot}:
\begin{align}
\nonumber
\E_x\left[\te{-q\t_a-r\ovl G_{\tau_a}}I_{A_o}\right]& =
\E_x\left[\te{-(q+r)T_{u,u+a}}I_{\left\{T_{u,u+a}=T^+_{u+a}\right\}}
\te{-q\left(\t_a-T_{u,u+a}\right)-r\left(\ovl G_{\t_a}-T_{u,u+a}\right)}I_{A_o}\right]\\
& = \frac{W^{(q+r)}(a\wedge (x-u))}{W^{(q+r)}(a)}\>
\E_{x\vee (u+a)}\left[\te{-q\t_a-r\ovl G_{\tau_a}}I_{A_o}\right].
\label{eq:First_Step}
\end{align}

\textit{Step (2)} By~\eqref{eq:First_Step}
we may assume that the starting point
$x$
of the process
$X$
is greater or equal to
$u+a$.
Since
$Y_{\tau_a-}\in[0,a]$
this implies that
\begin{equation}
\label{eq:Step_2_Observation}
\unl X_{\t_a-}\geq u\quad\P_x\text{-a.s. and hence}\quad
\left\{\unl X_{\t_a}\geq u\right\}=
\left\{X_{\t_a-}+\Delta X_{\tau_a}\geq u\right\}\>\P_x\text{-a.s.}
\end{equation}
In order to deal with the jump at time
$\tau_a$
we apply the compensation formula (see Section O.5 in~\cite{bert96})
to the Poisson point process
$(\Delta X_t)_{t\geq0}$
as follows:
\begin{align}
\nonumber
\E_x\left[\te{-q\t_a-r\ovl G_{\tau_a}}
I_{A_o}\right]& =
\E_x\left[\te{-q\t_a-r\ovl G_{\tau_a}}
I_{\left\{\ovl X_{\t_a}\in \td v, \unl X_{\t_a}\geq u, Y_{\tau_a-}\in\td y, Y_{\tau_a}-a\in\td h \right\}}\right]\\
\nonumber
& =
\E_x\left[\sum_{t\geq0}\te{-qt-r\ovl G_{t}}
I_{\left\{\sup_{s<t} \{Y_{s}\}\leq a, Y_{t}>a,
\ovl X_{t}\in \td v, Y_{t-}\in\td y\right\}}
I_{\left\{y-a-\Delta X_t\in\td h,  \Delta X_t\geq u+y-v\right\}}\right]\\
\nonumber & = I_{\{h\in(0,u-v-a)\}}\Lambda\left(y-a-\td h\right)
\E_x\left[\int_0^\infty \te{-qt-r\ovl G_{t}} I_{\left\{\ovl
X_{t}\in \td v, Y_{t-}\in\td y, \sup_{s<t} \{Y_{s}\}\leq
a\right\}}\td t
\right]\\
\label{eq:Final_Generalisation_of_resolvent} & =
I_{\{h\in(0,u-v-a)\}}\Lambda\left(y-a-\td h\right) (\E_x\times
E)\left[\te{-r\ovl G_{\eta}} I_{\left\{\ovl X_{\eta}\in \td v,
Y_{\eta}\in\td y, \eta<\tau_a\right\}} \right]/q,
\end{align}
where the sum runs over all $t\geq0$ such that $\Delta X_t\neq0$
and, as noted before, $q>0$, and $\eta$ is an exponential random
variable defined on a probability space $(\mbb R_+,\mc B(\mbb
R_+),P)$ such that $(\P\times P)[\eta>s|\mc F]=P(\eta>s)=\te{-qs}$
for $s\in\mbb R_+$. The second equality follows
by~\eqref{eq:Step_2_Observation}. The third equality follows from
the compensation formula for the Poisson point process of jumps of
$X$ and the fact that the inequality $\Delta X_t\geq u+y-v$ is
equivalent to the restriction on $h$ given in~\eqref{eq:ParmRest}.

\textit{Step (3)} The final task is to calculate the expectation
in~\eqref{eq:Final_Generalisation_of_resolvent}.
This requires the It\^o excursion theory for the reflected process
$Y$.
Since
$X$ is spectrally negative, the local time at zero of $Y$ can be
taken to be the running supremum of $X$, i.e. $L=\ovl X -x$, under
$\P_x$. Let $(t,\epsilon_t)_{t\geq0}$ be the Poisson point
process of excursions of
$X$
from its
supremum $\ovl X$, taking values in
$[0,\infty)\times(\EE\cup\{\partial\})$,
where
$$
\EE=\left\{\varepsilon\in D(\R):\exists \zeta\in(0,\infty]\text{ such that }\varepsilon(\zeta)=0 \text{ if } \zeta<\infty,\>
\varepsilon(0)\geq0,\>\varepsilon(t)>0\>\forall t\in(0,\zeta)\right\},
$$
$D(\R)$ is the Skorokhod space and $\partial$ is the graveyard
state, with the intensity measure $\td t\times n(\td
\varepsilon)$, where $n$ denotes the It\^o excursion measure. For
each instance of local time $t\in(0,v-x]$, the inverse local time
$L^{-1}$ satisfies $L_{L^{-1}_t}=\ovl X_{L^{-1}_t}-x=t$ since
$L^{-1}_t=T_{x+t}^+$. Since the excursions are indexed by local
time, for any $t\in(0,v-x]$ we have $\epsilon_t =
\{X_{L^{-1}_{t-}}- X_{L^{-1}_{t-}+s}\>:\>0<s\leq
L^{-1}_{t}-L^{-1}_{t-}\}$ if $\>L^{-1}_{t-}<L^{-1}_{t}$ and
$\epsilon_t=\partial$ if $\>L^{-1}_{t-}=L^{-1}_{t}$. This implies
that the excursions can be indexed by a subset of actual time that
is given by the left-end points of excursion intervals. For any
excursion $\varepsilon\in\EE$ define
$T_a(\varepsilon)=\inf\{s>0:\varepsilon(s)>a\}$ (with
the\label{page:Def_Of_T_a} convention $\inf\emptyset=\infty$) and
let $\zeta(\varepsilon)$ be the life time of $\varepsilon$. Define
also the height of the excursion $\varepsilon$ by
$\ovl\varepsilon=\sup\{\varepsilon(s):0<s< \zeta(\varepsilon)\}$.
We now obtain in the case $y>0$
\begin{eqnarray}
\nonumber \lefteqn{ (\E_x\times E)\left[\te{-r\ovl G_{\eta}}
I_{\left\{\ovl X_{\eta}\in \td v, Y_{\eta}\in\td y,
\eta<\tau_a\right\}}
\right]}\\
\nonumber
& = & (\E_x\times E)\left[\sum_{g\geq0}\te{-r\ovl G_{g}}
I_{\left\{\sup_{h<g}\left\{\ovl \epsilon_h\right\}\leq a,\>
\eta\in\left(g,g+(T_a\wedge\zeta)(\epsilon_g)\right),\>\ovl
X_g\in\td v\right\}} I_{\left\{\epsilon_g(\eta-g)\in\td y\right\}}
\right]\\
\nonumber
& = & (\E_x\times E)\left[\int_0^\infty \td \ovl X_s\te{-r\ovl
G_{s}} I_{\left\{\sup_{h<s}\left\{\ovl \epsilon_h\right\}\leq
a,\>\ovl X_s\in\td v\right\}}I_{\{\eta>s\}}\> n\left[t'<
(T_a\wedge\zeta)(\varepsilon),\>
\varepsilon(t')\in\td y\right]\big|_{t'=\eta-s}\right]\\
& = & \E_x\left[\int_0^\infty \te{-qs-r\ovl G_{s}}
I_{\left\{\sup_{h<s}\left\{\ovl \epsilon_h\right\}\leq a,\>\ovl
X_s\in\td v\right\}}\td \ovl X_s\right]  \nn \\ &\phantom{=}&
\cdot\ E\left[n\left(\eta< (T_a\wedge\zeta)(\varepsilon),\>
\varepsilon(\eta)\in\td y\right)\right],
\label{eq:Sum_In_Compensation_Formula}
\end{eqnarray}
where the sum in the first equality is over all left-end points of
excursion intervals. The second equality follows from the
compensation formula of excursion theory (see~\cite[Cor.
IV.11]{bert96}) and the third from Fubini's theorem, the
independence of $X$ and $\eta$ and the fact $P[\eta>s]=\te{-qs}$.
The second expectation in~\eqref{eq:Sum_In_Compensation_Formula}
is given by
\begin{eqnarray}
\label{eq:Stuff_From_Martijn} E\left[n\left(\eta<
(T_a\wedge\zeta)(\varepsilon),\> \varepsilon(\eta)\in\td
y\right)\right] = q\left[W^{(q)\prime}_+(y)-
\frac{W^{(q)\prime}_+(a)}{W^{(q)}(a)}W^{(q)}(y)\right]\td y \qquad
y\in(0,a).
\end{eqnarray}
The formula in~\eqref{eq:Stuff_From_Martijn} was proved
in~\cite{P} (see equation~(22)).

In order to obtain the first expectation in~\eqref{eq:Sum_In_Compensation_Formula}
first note that, since
at the time
$L^{-1}_t$
the process
$X$
is at its supremum level
$x+t$
(i.e. $L^{-1}_t=T_{x+t}^+$),
we must have
$\ovl G_{L^{-1}_t} = L^{-1}_t$
$\P_x$-a.s
for all
$t>0$.
Furthermore the identity
$\ovl X_{L^{-1}_t}-x=t$
implies that if we reparametrise the integral under the
first expectation in~\eqref{eq:Sum_In_Compensation_Formula}
using inverse local time
we find
\begin{eqnarray}
\nonumber
\lefteqn{
\E_x\left[\int_0^\infty \te{-qs-r\ovl G_{s}} I_{\left\{\sup_{h<s}\left\{\ovl \epsilon_h\right\}\leq
a,\>\ovl X_s\in\td v\right\}}\td \ovl X_s\right]}\\
\nonumber
& = & \E_x\left[\te{-(q+r)L^{-1}_{v-x} }
I_{\left\{\sup_{h<L^{-1}_{v-x} }\left\{\ovl \epsilon_h\right\}\leq a\right\}} \right]\,\td v\\
\nonumber
& = &
\te{-\Phi(q+r)(v-x)}
\E_x\left[\te{-(q+r)L^{-1}_{v-x}+\Phi(q+r)(v-x)}
I_{\left\{\sup_{h<L^{-1}_{v-x} }\left\{\ovl \epsilon_h\right\}\leq a\right\}} \right]\,\td v\\
& = &
\label{eq:Final_Prob}
\te{-\Phi(q+r)(v-x)}
\P_x^{\Phi(q+r)}\left[T_v^+<\tau_a\right]\,\td v.
\end{eqnarray}
The third equality follows by an application of the Esscher change
of measure formula and the equality in~\eqref{eq:Final_Prob} is a
consequence of the following
$$\{T_v^+<\tau_a\}= \{\forall t\in[0,v-x]\,\text{with}\,\epsilon_t\neq\partial,\,\ovl\epsilon_t\leq
a\}\quad\P_x\text{ -a.s.}$$ Since $N_t=\#\{u\in[0,t] :
\ovl\epsilon_u>a, \epsilon_u\neq\partial\}$ is a Poisson process
with parameter $n^{\Phi(q+r)}(\ovl\varepsilon>a)$ under the
measure $\P_x^{\Phi(q+r)}$, where $n^{\Phi(q+r)}$ is the It\^o
excursion measure of $X$ under the measure $\P_x^{\Phi(q+r)}$, we
obtain
\begin{eqnarray}
\label{eq:Final_Prob_Formula}
\P_x^{\Phi(q+r)}\left[T_v^+<\tau_a\right] & = & \te{-(v-x)n^{\Phi(q+r)}\left(\ovl \varepsilon>a\right)},
\qquad\text{where}\qquad\\
n^{\Phi(q+r)}\left(\ovl \varepsilon>a\right) & = &
\frac{W_{\Phi(q+r)+}'(a)}{W_{\Phi(q+r)}(a)}=
\frac{W^{(q+r)\prime}_+(a)}{W^{(q+r)}(a)}-\Phi(q+r)
\label{eq:PoissonRandomMEasureCompensatedFormulea_2}
\end{eqnarray}
and $W_{\Phi(q+r)}(a)$ denotes the 0-scale function under the
measure $\P_x^{\Phi(q+r)}$. The first equality
in~\eqref{eq:PoissonRandomMEasureCompensatedFormulea_2} is a
well-known representation of the scale function in terms of the
excursion measure and the second equality follows from the Esscher
change of measure formula
$W_{\Phi(q+r)}(x)=\te{-x\Phi(q+r)}W^{(q+r)}(x)$.

Since $\ovl G_\eta=\eta$ on the set $\{Y_\eta=0\}$ we find in the
case $y=0$
\begin{eqnarray*}
(\E_x\times E)[\te{-r G_\eta}I_{\{\ovl X_\eta>v, Y_\eta=0, \eta<
\tau_a\}}] &=& (\E_x\times E)[\te{-r \eta}I_{\{\ovl
X_\eta> v, Y_\eta=0, \eta< \tau_a\}}]\\
&=& \E_x\le[\int_0^\infty q \te{-(r+q) t} I_{\{T_v^+ < t, Y_t = 0,
t <
\tau_a\}}\td t \right]\\
&=& \frac{q}{r+q}\E_x[\te{-(r+q)T_v^+}I_{\{T_v^+<\tau_a\}}]
\E_0\le[\int_0^{\tau_a}(r+q)\te{-(q+r)t} I_{\{Y_t = 0\}}\td t\ri]\\
&=& q \exp\le(-v\frac{W^{(q+r)\prime}_+(a)}{W^{(q+r)}(a)}\ri)
\cdot \frac{W^{(q+r)}(a)}{W^{(q+r)\prime}_+(a)} W^{(q+r)}(0),
\end{eqnarray*}
where the first factor follows by combining~\eqref{eq:Final_Prob},
\eqref{eq:Final_Prob_Formula} and
\eqref{eq:PoissonRandomMEasureCompensatedFormulea_2} and the
second factor follows from~\cite[equation~(22)]{P}. Since
$W^{(q+r)}(0)=W^{(q)}(0)$ we thus obtain
\begin{equation}\label{eq:y=0}
(\E_x\times E)[\te{-r G_\eta}I_{\{\ovl X_\eta\in\td v, Y_\eta=0,
\eta< \tau_a\}}] = q
\exp\le(-v\frac{W^{(q+r)\prime}(a)}{W^{(q+r)}(a)}\ri) W^{(q)}(0).
\end{equation}

Identities~\eqref{eq:First_Step},~\eqref{eq:Final_Generalisation_of_resolvent},~\eqref{eq:Sum_In_Compensation_Formula},
\eqref{eq:Stuff_From_Martijn},~\eqref{eq:Final_Prob}
and~\eqref{eq:Final_Prob_Formula} and~\eqref{eq:y=0} together
imply the formula in~\eqref{eq:QuitLaq_overshoot}. This concludes
the proof in the case of the event $A_o$.

\smallskip
In the case of creeping,
step (1) consists of the application of the strong Markov property
and the first factor
in~\eqref{eq:QuitLaq_Creeping} follows from an analogous calculation
to~\eqref{eq:First_Step}
with
$A_o$
replaced by
$A_c$.
We therefore assume in what follows that
$x\geq u+a$.

Step~(2) does not feature in the context of $A_c$ as there is no
overshoot at $\tau_a$. For the analogue of step~(3) note that
$x\geq u+a$ implies (see also Figure~\eqref{fig:Main}) the
inclusion
$$
\left\{Y_{\tau_a}=a\right\}\subset
\left\{\unl X_{\tau_a}\geq u\right\}\quad\P_x\text{ -a.s.}
$$
The following excursion calculation, similar to~\eqref{eq:Sum_In_Compensation_Formula},
provides a key step:
\begin{eqnarray}
\nonumber
\lefteqn{\E_x\left[\te{-q\t_a-r\ovl G_{\tau_a}}I_{\left\{\ovl X_{\tau_a}\in \td v,
Y_{\tau_a}=a\right\}}\right]}\\
& = &
\label{eq:Second_Creeping}
\E_x\left[\sum_{g\geq0}\te{-gq-r\ovl G_g} I_{\left\{\sup_{h<g}\left\{\ovl
\epsilon_h\right\}\leq a,\>\ovl X_{\tau_a}\in \td v \right\}}\te{-qT_a(\epsilon_g)}
I_{\left\{T_a(\epsilon_g)<\zeta(\epsilon_g),\epsilon_g(T_a(\epsilon_g))=a\right\}}\right] \\
& = &
\E_x\left[\int_0^\infty \te{-qs-r\ovl G_{s}} I_{\left\{\sup_{h<s}\left\{\ovl \epsilon_h\right\}\leq
a,\>\ovl X_s\in\td v\right\}}\td \ovl X_s\right]\,
\int_{\EE} \te{-qT_a(\varepsilon)} I_{\left\{T_a(\varepsilon)<\zeta(\varepsilon),\varepsilon(T_a(\varepsilon))=a\right\}}\,n(\td \varepsilon)
\label{eq:Final_Factorization_Creeping}
\end{eqnarray}
The summation in~\eqref{eq:Second_Creeping}
is over the left-end points
of the excursion intervals and the equality in~\eqref{eq:Final_Factorization_Creeping}
follows from the compensation formula of It\^o's excursion theory.

The expectation in the first factor
in~\eqref{eq:Final_Factorization_Creeping} is
identical to~\eqref{eq:Final_Prob}.
To find the second factor in~\eqref{eq:Final_Factorization_Creeping}
first note for any
$x>0$
the following identities hold as a consequence of the
compensation formula:
\begin{eqnarray*}
\lefteqn{\E\left[\te{-qT^-_{-a}} I_{\left\{X_{T^-_{-a}}=-a, T^-_{-a}<T^+_x\right\}}\right]}\\
& = &
\E\left[\sum_{g\geq0}\te{-gq} I_{\left\{\sup_{h<g}\left\{\ovl \epsilon_h-\ovl X_{h}\right\}\leq a,\>\ovl X_{g}\leq x \right\}}
\te{-qT_{a+\ovl X_{g}}(\epsilon_g)}
I_{\left\{T_{a+\ovl X_{g}}(\epsilon_g)<\zeta(\epsilon_g),\epsilon_g(T_{a+\ovl X_{g}}(\epsilon_g))=a+\ovl X_{g}\right\}}\right] \\
& = &
\E\left[\int_0^x \te{-qL^{-1}_y}
I_{\left\{\sup_{h<L^{-1}_y}\left\{\ovl \epsilon_h-\ovl X_h\right\}\leq a \right\}}
\td y\,
\int_{\EE}
\te{-qT_{a+y}(\varepsilon)}
I_{\left\{T_{a+y}(\varepsilon)<\zeta(\varepsilon),\varepsilon(T_{a+y}(\varepsilon))=a+y\right\}}n(\td
\varepsilon)\right].
\end{eqnarray*}
The second equality follows by reparametrising the integral
which arises in the compensation formula
using inverse local time.
Since
$L^{-1}_0=0$ a.s.
and the trajectories of
$X$
and
$L^{-1}$
are right-continuous, the following equality holds:
\begin{eqnarray}
\label{eq:Final_Creepy}
\int_{\EE} \te{-qT_a(\varepsilon)}
I_{\left\{T_a(\varepsilon)<\zeta(\varepsilon),
\varepsilon(T_a(\varepsilon))=a\right\}}\,n(\td \varepsilon)
&=& \lim_{x\searrow0} x^{-1}
\E\left[\te{-qT^-_{-a}} I_{\left\{X_{T^-_{-a}}=-a, T^-_{-a}<T^+_x\right\}}\right].
\end{eqnarray}
This limit can also be derived from Proposition 2 in \cite{Doney}.
Insert identity~\eqref{eq:twosmoothdown}
into~\eqref{eq:Final_Creepy}, apply L'Hopital's rule and the fact
$W^{(q)\prime}(0)=\frac{2}{\sigma^2}$ to obtain a formula for the
second factor in~\eqref{eq:Final_Factorization_Creeping}. This
implies~\eqref{eq:QuitLaq_Creeping} and therefore concludes the
proof of the theorem. \exit

\subsection{Proof of Theorem~\ref{thm:Triple}}
We start by noting that the spectral negativity of
$X$
implies the following identity
$X_{\WH \t_a}- \unl X_{\WH \t_a}=\WH Y_{\WH \t_a} = a$
a.s.
We can therefore apply the Esscher transform to~\eqref{eq:Triple_Law_Equation}
to obtain
\begin{equation}
\label{eq:Triple_Law_Equation_MAesure_Changed} \E\left[
\te{-q\WH \t_a - r\unl G_{\WH\tau_a}+ u\unl X_{\WH \t_a}}
I_{\left\{\ovl X_{\WH \t_a}< v\right\}} \right] = \te{-au}
\E^u\left[ \te{-(q-\psi(u))\WH \t_a - r\unl G_{\WH\tau_a}}
I_{\left\{\ovl X_{\WH \t_a}< v\right\}} \right].
\end{equation}
If $a<v$, then $\P\left[\ovl X_{\WH \t_a}< v\right]=1$ as
remarked above (see also Figure~\ref{fig:Min_Fig}).
Note further that
$$\left\{\ovl X_{\WH \t_a}<v\right\}=\left\{T_{v-a}^-<T_v^+\right\}\qquad\text{for any}\quad v\in[0,a].$$
Observe
that $X_{T_{v-a}^-}\leq v-a$ $\P$-a.s.
Hence the strong Markov property
applied to the right-hand side
of~\eqref{eq:Triple_Law_Equation_MAesure_Changed} and
formula~\eqref{eq:two-sided2} yield
\begin{align*}
\E^u\left[ \te{-(q-\psi(u))\WH \t_a - r\unl G_{\WH\tau_a}}
I_{\left\{\ovl X_{\WH \t_a}< v\right\}} \right]
& =
\E^u\left[ \te{-(q+r-\psi(u))T_{v-a}^-} I_{\left\{T_{v-a}^-<T_v^+\right\}}\right]
\E^u\left[\te{-(q-\psi(u))\WH \t_a - r\unl G_{\WH\tau_a}}\right] \\
& = 
\left[ Z_u^{(p)}(a-v)- Z_u^{(p)}(a)
\frac{W_u^{(p)}(a-v)}{W_u^{(p)}(a)} \right]
\frac{W^{(q+r)}(a)}{W^{(q)}(a)Z^{(p)}_u(a)}
\end{align*}
for all
$v\geq0$,
where $p=q+r-\psi(u)$. To complete the
proof we need to show that
\begin{equation}\label{eq:tg_under_u}
\E^u\left[\te{-(q-\psi(u))\WH \t_a - r\unl G_{\WH\tau_a}}\right] =
\frac{W^{(q+r)}(a)}{W^{(q)}(a)Z^{(p)}_u(a)}
\quad\text{for}\quad q,r,u\geq0.
\end{equation}
This, together
with~\eqref{eq:scala_Functions_W}--\eqref{eq:scala_Functions_Z},
implies the formula in the theorem.

By the compensation formula applied to the Poisson point process
of excursions of $\WH Y$ away from zero,
\begin{equation}
\label{eq:Lucky_Identity}
\E^u[\te{-(q-\psi(u))\WH\tau_a - r\unl G_{\WH\tau_a}}]  =  \E^u\le[\int_0^\infty
\te{-pt}I_{\{t<\unl G_{\WH\tau_a}\}}\td \WH L_t\ri] \int_\EE
\te{-(q-\psi(u)) T_a(\varepsilon)}I_{\{T_a(\varepsilon)<\zeta(\varepsilon)\}} \WH n^u(\td\varepsilon),
\end{equation}
where $\WH n^u$ is the corresponding excursion measure under the
probability measure $\P^u$, $\WH L$ a local time process of $\WH
Y$ at zero and $T_a$ and $\zeta$ are as in the proof of
Theorem~\ref{thm:lawreflect1} (see
page~\pageref{page:Def_Of_T_a}). From Bertoin~\cite[Ch.~VII,
Prop.~15]{bert96} it follows that for any $(\mc F_t)$-stopping
time $\tau$ and $A\in\mc F_\tau$ we have
\begin{equation}
\label{eq:Final_Hopefully}
\WH n^u(A \cap\{\t<\zeta\})=\lim_{x\searrow 0}W_u(x)^{-1}\E^u_x\le[I_{A\cap\{\t<T_0^-\}}\ri]
\qquad\text{for any}\quad u\geq0.
\end{equation}
This limit can also be derived from \cite[Proposition 2]{Doney}.
We now apply~\eqref{eq:Final_Hopefully} to the stopping time
$\tau=T_a^+$ to obtain
\begin{eqnarray}
\nonumber
\int_\EE \te{-(q-\psi(u)) T_a(\varepsilon)}I_{\{T_a(\varepsilon)<\zeta(\varepsilon)\}} \WH n^u(\td\varepsilon) &=&
\lim_{x\searrow 0}
W_u(x)^{-1}\E^u_x\le[\te{-(q-\psi(u)) T_a^+}I_{\{T_a^+<T_0^-\}}\ri]\\
\nonumber
&=& W_u^{(q-\psi(u))}(a)^{-1} \lim_{x\searrow 0} \frac{W_u^{(q-\psi(u))}(x)}{W_u(x)}\\
& =&
W_u^{(q-\psi(u))}(a)^{-1}
\label{eq:Compensation_Integral_Final}
\end{eqnarray}
for any
$q\geq\psi(u)$.
The second equality follows from~\eqref{eq:two-sided1}.
The third equality is a consequence of the well-known identity
$W^{(l)}(x)=\sum_{k=0}^\infty l^k W^{*(k+1)}(x)$,
which holds for all
$l\in\C$,
$x\geq0$,
where
$W^{*k}$
denotes the
$k$-th convolution of
$W$
with itself.

The identity in~\eqref{eq:Lucky_Identity}
for
$r=0$,
the Laplace transform
$\E^u\left[ \te{-(q-\psi(u))\WH \t_a} \right] = 1/Z_u^{(q-\psi(u))}(a)$,
given in~\cite[Prop.~2]{P},
and
identity~\eqref{eq:Compensation_Integral_Final}
imply
$$
\E^u\le[\int_0^\infty \te{-qt}I_{\{t<\unl G_{\WH\tau_a}\}}\td \WH L_t\ri]
= \frac{W_u^{(q-\psi(u))}(a)}{Z_u^{(q-\psi(u))}(a)}
$$
for
$q\geq\psi(u)$.
This identity, together with~\eqref{eq:Lucky_Identity}
and~\eqref{eq:Compensation_Integral_Final},
yields~\eqref{eq:tg_under_u}
for large
$q$.
Analyticity of expressions in~\eqref{eq:tg_under_u}
on both sides of the identity imply~\eqref{eq:tg_under_u}
for all parameter values. This concludes the proof of the
theorem.
\exit

\end{document}